\documentclass[11pt]{amsart}
\usepackage{amsmath}
\usepackage{amsfonts}
\usepackage{latexsym}
\usepackage{graphicx}

\newcommand \R {{ \mathbb R}}
\newcommand \C {{ \mathbb C}}
\newcommand \Z {{ \mathbb Z}}

\newtheorem{theorem}{Theorem}[section]
\newtheorem {lemma} [theorem]{Lemma}
\newtheorem {proposition}[theorem]{Proposition}

\newtheorem{claim}[theorem]{Claim}

\title[Cohomological Equation for Horocycle Maps]
{The Cohomological Equation and Invariant distributions for Horocycle Maps}

\author{James Tanis} 

\address{Department  of Mathematics\\
  University of Maryland \\
  College Park, MD USA}

\begin{document}

  \begin{abstract}
We study the invariant distributions for horocycle maps on $\Gamma\backslash SL(2, \mathbb{R})$ and prove Sobolev estimates for the cohomological equation of horocycle maps.  As an application, we obtain a rate of equidistribution for horocycle maps on compact manifolds.  
  \end{abstract}
  \maketitle

\tableofcontents

\section{Introduction}

We say that $f$ is a $\textit{coboundary for the flow}$ $\{\phi_t\}_{t\in \mathbb{R}}$ if there is a solution $u$ to the $\textit{cohomological}$ $\textit{equation}$ $$\frac{d}{dt} u\circ \phi_t |_{_{t = 0}} = f,$$ and $f$ is a$\textit{ coboundary for the map}$ $\phi_{T}$ if there is a solution to the $\textit{cohomological}$ $\textit{equation}$ $$u\circ \phi_{T} - u = f,$$ for $T > 0$.  In this paper we study the discrete analogue of the classical horocycle flow $\Big\{\left(\begin{array}{rr} 1 & t\\ 0 &1\end{array}\right)\Big\}_{t \in \mathbb{R}}$ called classical horocycle maps $\Big\{\left(\begin{array}{rr} 1 & n T\\ 0 &1\end{array}\right)\Big\}_{n \in \mathbb{Z}},$ each acting by right multiplication on (compact) homogeneous spaces of the form $\Gamma\backslash PSL(2, \mathbb{R})$.  
Motivated by the success of using cohomological equations to prove quantitative equidistribution of horocycle flows and nilflows \cite{2}, \cite{3}, we study the cohomological equation for horocycle maps and quantitative equidistribution. 

 Horocycle flows are known to have zero entropy, and the precise mixing rates for geodesic and horocycle flows were obtained by Ratner \cite{16} and Moore \cite{17}, and Ratner proved horocycle flows have polynomial decay of correlations.  Concerning ergodicity, Furstenberg \cite{15} proved the horocycle flow is uniquely ergodic (i.e. every orbit equidistributes) in 1970.  M. Burger \cite{30} estimated the rate of unique ergodicity for sufficiently smooth functions along orbits of horocycle flows on compact surfaces and on open complete surfaces of positive injectivity radius.  P. Sarnak \cite{24} obtained asymptotics for the rate of unique ergodicity of cuspidal horocycles on noncompact surfaces of finite area using a method based on Eisenstein series.  For sufficiently regular functions, Flaminio-Forni \cite{2} improved on Burger's estimate for compact surfaces by establishing precise asymptotics in this setting, and in the case of noncompact, finite area surfaces, they generalize the result of P. Sarnak to arbitrary horocycle arcs.
 
 Quantitative equidistribution results for horocycle maps are very recent.  Shah's conjecture states that for all $\delta > 0$, the horocycle map $\{\phi_{n^\delta}^U\}_{n\in \mathbb{Z}^+}$ equidistributes in $\Gamma \backslash SL(2, \mathbb{R})$.  In \cite{23}, Venkatesh was able to use quantitative equidistribution and quantitative mixing of the horocycle flow to prove upper bounds on the equidistribution rate of the "twisted" horocycle flow $\{\phi_t^U\times e^{2\pi i t}\}_{t\in \mathbb{R}}$ on compact manifolds $\Gamma\backslash SL(2, \R) \times S^1$.  He then used this to estimate a rate of equidistribution for the horocycle map $\{\phi_{n}^U\}_{n \in \mathbb{Z}^+}$ in compact $\Gamma\backslash SL(2, \mathbb{R})$ and prove $\{\phi_{n^{(1 + \delta)}}^U\}_{n \in \mathbb{Z}^+}$ equidistributes whenever $0 \leq \delta < \delta_{\Gamma}$, for some explicit number $0 < \delta_{\Gamma} << 1.$  Even more recently, Sarnak-Ubis proved the orbit of a generic point at prime times of the horocycle flow in the modular surface is dense in a set of positive measure \cite{6}.
 
Representation theory is a natural tool for studying cohomological equations on homogeneous spaces \cite{2}, \cite{3}, \cite{10}
.  Flaminio-Forni's (2003, \cite{2}) detailed analysis of the cohomological equation for the horocycle flow was carried out through its representations in the irreducible, unitary components of $L^2(\Gamma \backslash SL(2, \mathbb{R}))$.  We take this approach for the cohomological equation of horocycle maps.  
Previous results on cohomological equations for homogeneous $\mathbb{R}$ or $\mathbb{Z}$ actions show there are infinitely many independent distributional obstructions to the existence of $L^2$ solutions.
Consistent with this picture, we find there are infinitely many independent distributional obstructions for the horocycle map with some finite loss of regularity between the Sobolev estimates of the transfer function and its coboundary (see also \cite{2}, \cite{3}, \cite{9}), and we find this for horocycle maps as well.

We obtain an asymptotic formula for the ergodic sum of the horocycle map in terms of the invariant distributions, and we improve the estimate for the rate of equidistribution found in \cite{23} for compact manifolds.  As in \cite{2}, we use our estimate of the cohomological equation for the map to obtain a rate of equidistribution for coboundaries, and we use the analysis of the flow invariant distributions for the horocycle flow in \cite{2} to estimate the rate of decay for the flow invariant distributions of the map.  We use Venkatesh's estimate of the equidistribution of the twisted horocycle flow in \cite{23} to estimate the invariant distributions of the map that are not flow invariant.  Then because the ergodic sum of every regular enough function is controlled either by the cohomological equation or one of the invariant distributions, we obtain an upper bound on the speed of equidistribution. 

\subsection{Preliminary definitions}

The Poincar$\acute{\text{e}}$ upper half-plane $H$ is the manifold $\{z \in \C  |\Im(z) > 0\}$ endowed  with the metric $\frac{|dz|^2}{(\Im z)^2}.$  Its isometry group is $PSL(2, \R)$.  If $\Gamma \subset PSL(2, \mathbb{R})$ is a discrete subgroup acting without fixed points, then $M := \Gamma\backslash H$ is a Riemannian manifold of constant negative curvature.  Let $SH$ be the unit tangent bundle of $H$.  As $PSL(2, \R)$ acts simply transitively on $SH$, we obtain the identifications $PSL(2, \R) \approx SH,$ and $\Gamma\backslash SH \approx \Gamma\backslash PSL(2, \mathbb{R}).$  Define $SM = \Gamma\backslash SL(2, \R)$, which is therefore the double cover of $\Gamma\backslash SH$.

The matrices \begin{equation}\label{370}U = \left(\begin{array}{rr}0 & 1\\ 0 & 0\end{array}\right) \text{ and } V = \left(\begin{array}{rr}0 & 0\\ 1 & 0\end{array}\right),\end{equation} in $sl(2, \mathbb{R})$ are the stable and unstable "horocycle vector fields" on $SM$; the right multiplication on $\Gamma\backslash SL(2, \R)$ by the one-parameter groups $$\phi_t^U := e^{t U} = \left(\begin{array}{rr}1 & t\\0 & 1\end{array}\right) \text{ and } \phi_t^V := e^{t V} = \left(\begin{array}{rr}1 & 0\\t & 1\end{array}\right)$$ defines the flows $\{\phi_t^U\}_{t\in \R}$ and $\{\phi_t^V\}_{t\in \R}$ that correspond to the stable and unstable horocycle flows on $SM$, respectively.  Let $T > 0$, and define $$L_T u := u \circ \phi_T^U - u,$$ for $u\in L^2(\R)$.  The main result in this paper is to study the cohomological equation \begin{equation}\label{716}L_T u = f\end{equation} when $f$ and $u$ belong to some Sobolev spaces to be made precise later.  In particular, we are interested in obtaining Sobolev estimates of the transfer function $u$ in terms of the coboundary $f$.
\subsection{Harmonic analysis}

Elements of $sl(2, \mathbb{R})$ generate some area preserving flows on $SM$, and we choose a basis for $sl(2, \mathbb{R})$ to be 
\begin{equation}\label{371}X = \left(\begin{array}{rr}1 & 0 \\0 & -1\end{array}\right), \ \ Y = \left(\begin{array}{rr}0 & -1 \\-1 & 0\end{array}\right), \ \ \Theta = \left(\begin{array}{rr}0 & 1 \\-1 & 0\end{array}\right),\end{equation} which are generators for the geodesic, orthogonal geodesic and circle vector fields respectively.  These generators satisfy the commutation rules $$[X, Y] = -2 \Theta, \ \ \ [\Theta, X] = 2 Y,  \ \ \ [\Theta, Y] = -2 X,$$ and note that we have $$U = \frac{\Theta - Y}{2} \text{ and } V = \frac{-Y - \Theta}{2}.$$  

We define the$\textit{ Laplacian }\triangle$ as the element of the enveloping algebra of $sl(2, \mathbb{R})$ given by $$\triangle := - X^2 - Y^2 - \Theta^2.$$  The$\textit{ Casimir }$operator is then given by $$\Box := \triangle + 2 \Theta^2;$$ it generates the center of the enveloping algebra of $sl(2, \mathbb{R}).$  As such, it acts as a constant $\mu \in \mathbb{R}$ on each irreducible, unitary representation space $\mathcal{K}_\mu$, 
and its value classifies the $\mathcal{K}_\mu$ into three classes.  The representation $\mathcal{K}_\mu$ belongs to the $\textit{principal}$ $\textit{series}$ if $\mu\geq 1$, the$\textit{ complementary series }$if $0 < \mu < 1$, the$\textit{ discrete series }$if $\mu\in \{-4j^2 + 4j| j \geq 2 \text{ is an integer}\}$ and the$\textit{ mock discrete series }$if $\mu = 0.$  

Our notation differs from other conventions in the following two ways.  We consider the discrete series and the mock discrete series together, so for any $\mu \leq 0$, we simply refer to $\mathcal{K}_\mu$ as a discrete series component.  Next, our measure for the holomorphic unit disk model of the discrete series is $4 \frac{(1 - |\xi|^2)^{\nu - 1}}{|\xi - 1|^{2(\nu + 1)}} du dv,$ which is obtained by change of variable from the corresponding measure $y^{\nu - 1} dx dy$ in the upper half-plane model.   Finally, some authors scale the vector fields so that the geodesic flow travels at unit speed with respect to the hyperbolic metric of constant curvature -1, and in this case, the component $\mathcal{K}_\mu$ is in the principal series whenever $\mu \geq 1/4$, e.g. \cite{2}.  Our geodesic flow has speed 2 with respect to the hyperbolic metric of constant curvature -1.

Considering $L^2(M)$ as a subspace of $L^2(SM)$, we have that $L^2(M)$ is $\Box$-invariant, and the spectrum of $\Box$ on $L^2(M)$ coincides with the spectrum of the Laplace-Beltrami operator on $\triangle_M$.  The spectrum of the Laplace-Beltrami operator $\triangle_M$ on $M$ and that of the Casimir $\Box$ coincide on $\mathbb{R}^+$.  When $M$ is compact, standard elliptic theory shows $spec(\triangle_M)$ is pure point and discrete, with eigenvalues of finite multiplicity.  When $M$ is not compact, $spec(\triangle_M)$ is Lebesgue on $[1, \infty)$ with multiplicity equal to the number of cusps, has possibly embedded eigenvalues of finite multiplicity in $[1, \infty)$, and has at most finitely many eigenvalues of finite multiplicity in (0, 1) (see \cite{24}).

There is a standard unitary representation of $SL(2, \mathbb{R})$ on the separable Hilbert space $L^2(SM)$ of square integrable functions with respect to the $SL(2, \mathbb{R})$ invariant volume form on $SM$.  As in Flaminio-Forni (2003), the Laplacian gives unitary representation spaces a natural Sobolev structure.  The $\textit{ Sobolev space of}\textit{ order }$ $\textit{r }$ $> 0$ is the Hilbert space $W^r(SM)\subset L^2(SM)$ that is the maximal domain determined by the inner product $$\langle f, g\rangle_{W^r(SM)} := \langle (1 + \triangle)^r f, g\rangle_{L^2(M)}$$ for $f, g\in L^2(SM)$.

The space of$\textit{ infinitely differentiable functions }$is $$C^{\infty}(SM) := \cap_{r \geq 0} W^r(SM).$$  For $r > 0$, the distributional dual to $W^r(SM)$ is the Sobolev space $W^{-r}(SM) = \left(W^{r}(SM)\right)'.$  The distributional dual to $C^{\infty}(SM)$ is $$\mathcal{E}'(SM) := \left(C^{\infty}(SM)\right)'.$$ 

Because the Casimir operator is the center of the enveloping algebra and acts as an essentially self-adjoint operator, any non-trivial unitary representation $\mathcal{H}$ for $SL(2, \mathbb{R})$ has a $SL(2, \mathbb{R})$-invariant direct integral decomposition \begin{equation}\label{intdec}\mathcal{H} = \int_{\mu \in spec(\Box)} \mathcal{K}_\mu d\beta(\mu),\end{equation} where $d\beta(\mu)$ is a Stiltjes measure over the spectrum $spec(\Box)$ (see\cite{25}). The space $\mathcal{K}_\mu$ does not need to be irreducible but is generally a direct sum of an at most countable number of equivalent irreducible components given by the spectral multiplicity of $\mu \in spec(\Box).$  

Additionally, all operators in the enveloping algebra are decomposable with respect to the direct integral decomposition ($\ref{intdec}$). In particular, $$L^2(SM) = \int_{\oplus_\mu} \mathcal{K}_\mu,$$ and for $r\in \mathbb{R}$, \begin{equation}\label{304}W^r(SM) = \int_{\oplus_\mu} W^r(\mathcal{K}_\mu).\end{equation}  When $SM$ is compact, these integral decompositions are direct sums.

\subsection{Statement of results}
$$ $$

Let $$\mu_0 = \inf(spec(\triangle_M) -\{0\}).$$  We consider manifolds $M$ with a $\textit{spectral gap}$ $\mu_0 > 0$.  Let $(x, T)\in SM\times \mathbb{R}^+$ and $r > 0$.  Let $$\mathcal{I}(SM) := \{\mathcal{D}\in \mathcal{E}'(SM): L_{-T} \mathcal{D} = 0\}$$ be the space of $\phi_T^U$-invariant distributions, and let $\mathcal{I}(\mathcal{K}_\mu) := \mathcal{I}(SM)\cap \mathcal{E}'(\mathcal{K}_\mu).$  Similarly, $$\mathcal{I}^r(SM) := \mathcal{I}(SM)\cap W^{-r}(SM) \text{ and } \mathcal{I}^r(\mathcal{K}_\mu) := \mathcal{I}(\mathcal{K}_\mu)\cap W^{-r}(\mathcal{K}_\mu).$$  By $(\ref{304})$, we have $$\mathcal{I}(SM) = \int_{\oplus} \mathcal{I}(\mathcal{K}_\mu) \text{ and } \mathcal{I}^r(SM) = \int_{\oplus} \mathcal{I}^r(\mathcal{K}_\mu).$$

\begin{theorem}\label{322}
Let $\epsilon > 0$ and $\mu_0 > 0$.  For all $\mu\in spec(\Box)$, the space $\mathcal{I}(\mathcal{K}_\mu)$ has infinite countable dimension.   

For $\mu > 0$, $\mathcal{I}(\mathcal{K}_\mu) \subset W^{-((1 + \Re\sqrt{1 - \mu})/2 + \epsilon)}(\mathcal{K}_\mu).$ 

When $\mu \leq 0,$ there is an infinite basis $\{\mathcal{D}_k\}_{k\in \mathbb{N}} \cup \{\mathcal{D}^0\} \subset \mathcal{I}(\mathcal{K}_\mu)$ such that $\mathcal{D}^0\in W^{-((1 + \Re\sqrt{1 - \mu})/2 + \epsilon)}(\mathcal{K}_\mu)$ is the flow invariant distribution studied in \cite{2} and $\langle \{\mathcal{D}_k\}_{k\in \mathbb{N}}\rangle \subset W^{-(1 + \epsilon)}(\mathcal{K}_\mu)$.
\end{theorem}

It will follow from Theorem $\ref{93}$ that the invariant distributions classify the space of coboundaries that have smooth solutions.  Let $$Ann(\mathcal{I}^r(SM)) = \{ f\in W^r(SM) | D f = 0\text{ for all }D\in \mathcal{I}^r(SM)\}.$$

\begin{theorem}\label{93}
Let $T > 0, r \geq 0$ and $f\in Ann(\mathcal{I}^{3r + 4}(SM))$.  Then there is a unique $L^2(SM)$ solution $u$ to $L_T u = f,$ and there is a constant $C_{r, T, SM} > 0$ such that \begin{equation}\label{301}\| u \|_{W^r(SM)} \leq C_{r, T, SM}
\| f \|_{W^{3r + 4}(SM)}.\end{equation}
\end{theorem}

If $\mathcal{D}$ is an invariant distribution and $u \in C^{\infty}(SM)$, then from definitions we conclude $$\mathcal{D}(f) = \mathcal{D}(u\circ \phi_T^U) - \mathcal{D}(u) = 0.$$  In this sense, invariant distributions are $\textit{ obstructions }$ to the existence of smooth solutions of equation ($\ref{716}$).  Theorem $\ref{302}$ gives the invariant distributions that obstruct the existence of $L^2(SM)$ solutions for sufficiently regular coboundaries $f$.  Let $$\tilde{\mathcal{I}}(SM) := \mathcal{I}(SM) - \int_{\oplus_\{\mu < 0\}} \langle \{\mathcal{D}_\mu^0\}\rangle d\beta(\mu).$$

\begin{theorem}\label{302}
Let $\mu_0 > 0$, $f\in W^9(SM)$ and $\mathcal{D}\in \tilde{\mathcal{I}}(SM)$.  If there exists $u \in L^2(SM)$ such that $L_T u = f$, then $\mathcal{D}(f) = 0.$  

Moreover, if $r > 1$ and $\mathcal{D} \in \int_{\oplus_{\{\mu < 0\}}} \langle \{\mathcal{D}_\mu^0\}\rangle d\beta(\mu),$ then there exists $f \in W^r(SM)$ and $u \in L^2(SM)$ such that $\mathcal{D}(f) \neq 0$ and $L_T u = f$.
\end{theorem}

In light of Theorem $\ref{93}$, Theorem $\ref{302}$ says, for example, that $L^2(SM)$ solutions for $C^{\infty}(SM)$ coboundaries with no discrete series component are automatically $C^{\infty}(SM)$.

We prove estimate $(\ref{301})$ on every irreducible component and then glue the solutions together.  Explicitly, suppose we are given $0 \leq r < t$, $\{u_\mu\}_\mu, \{f_\mu\}_\mu \in \int_{\oplus_\mu} \mathcal{K}_\mu$ and a constant $C_{r, t} > 0$ such that for all $\mu\in spec(\Box)$, \begin{equation}\label{137}\| u_{\mu}\|_{W^r(\mathcal{K}_\mu)} \leq C_{r, t} \| f_\mu \|_{W^{t}(\mathcal{K}_{\mu})}.\end{equation} Write $$u = \int_{\oplus_\mu} u_\mu, \ \ \ \ \ f = \int_{\oplus_\mu} f_\mu,$$ and observe $$\|u\|_{W^r(SM)}^2 = \|\int_{\oplus_\mu} u_{\mu}\|_{W^r(\mathcal{K}_\mu)}^2 = \int_{\oplus_\mu} \|u_{\mu}\|_{W^r(\mathcal{K}_\mu)}^2$$ $$ \leq C_{r, t}^2 \int_{\oplus_\mu} \|f_{\mu}\|_{W^{t}(\mathcal{K}_{\mu})}^2 = C_{r, t}^2 \| f \|_{W^{t} (SM)}^2.$$ 

It therefore suffices to establish $(\ref{137})$.  

The$\textit{ key idea }$to obtain estimate ($\ref{137})$ is to introduce a finite dimensional space $Y$ of additional distributions with the property that whenever a function is in $Ann(Y)$, the estimate $(\ref{137})$ is substantially easier to prove.  Then we remove these distributions using a dual basis to $Y$ consisting of explicit coboundaries and obtain $(\ref{137})$ for each dual basis element.  Combining gives the full estimate.\\

Horocycle maps and the horocycle flow are related through the following proposition:

\begin{proposition}\label{333}
Let $\triangle_{M}$ have a spectral gap, 
 let $s > 1$, $f\in W^{s}(SM)$ and \begin{equation}\label{320}A_T(f) := \int_0^T f\circ \phi_t^U dt.\end{equation}  Then there exists $u\in L^2(SM)$ such that $$\mathcal{L}_U u = f \text{ if and only if } u \circ \phi_T^U - u = A_T f.$$
\end{proposition} 
One can prove (using Theorem 1.2, Proposition $\ref{333}$ and Lemma $\ref{2070}$ together with Theorem 1.2 of \cite{2}) that the operator $A_T$ maps the space of smooth coboundaries for the flow bijectively onto the space of smooth coboundaries for the horocycle map.  Thus, one could possibly reduce the study of the cohomological equation for the time-$T$ map to that of the flow studied in Flaminio-Forni \cite{2}, where now obtaining the Sobolev estimate ($\ref{301}$) is equivalent to proving a lower Sobolev bound on the operator $A_T$.   However, this approach does not seem any easier, so we study time-$T$ maps directly.  In this way our results on the cohomological equation are completely independent of Flaminio-Forni. \\ \\

As an application to the above analysis, we prove a rate of equidistribution for horocycle maps.  Let $\alpha(\mu_0) = \frac{(1 - \sqrt{1 - \mu_0})^2}{8(3 - \sqrt{1 - \mu_0})},$ where $\mu_0 > 0$ is the spectral gap.  Define $\mathcal B_\mu := \{\mathcal D_k\}_{k} \cup \{\mathcal D^0\} \subset \mathcal{I}(\mathcal K_\mu)$ and $\mathcal B_\mu^s := \mathcal B_\mu \cap \mathcal I^s(\mathcal K_\mu).$  For all $\mu \in spec(\Box)$ and $\mathcal{D} \in \{\mathcal{D}_k\}_k\cup \{\mathcal{D}^0\}\subset \mathcal{I}(\mathcal{K}_\mu)$, define $$\mathcal{S}_\mathcal{D} := \left\{\begin{array}{llll}\alpha(\mu_0) \text{ if } \mathcal{D} = \mathcal{D}_k, \ k \neq 0\\ \frac{1 - \Re\sqrt{1 - \mu}}{2} \text{ if } \mathcal{D} = \mathcal{D}_0 \\ \frac{1 + \Re\sqrt{1 - \mu}}{2} \text{ if } \mathcal{D} = \mathcal{D}^0 \text{ and } \mu > 0 \\ 1 \text{ if } \mathcal{D} = \mathcal{D}^0 \text{ and } \mu \leq 0.\end{array}\right.$$   

\begin{theorem}\label{152}
Let $\phi_1^U$ be the horocycle map on the unit tangent bundle $SM$ of a compact hyperbolic Riemann surface $M$ with spectral gap $\mu_0 > 0$, and let $s \geq 14$.  Then there is a constant $C_{s} > 0$ such that for all $f\in W^{s}(SM)$ with zero average and $(x_0, N)\in SM\times \mathbb{Z}^+$, we have $$\frac{1}{N}\sum_{k = 0}^{N - 1} f(\phi_k^Ux_0)$$ \begin{equation}\label{155} = \sum_{\mu \in spec(\Box)} \sum_{\mathcal{D}\in \mathcal B_\mu^s} c_\mathcal{D}(x_0, N, s)\mathcal{D}(f) N^{-\mathcal{S}_{\mathcal{D}}} \log^+(N) + \mathcal{R}(x, N, s)(f),\end{equation} where the remainder distribution $\mathcal{R}(x_0, N, s)$ is an element of $\cap_{\mu \in spec(\Box)}(\mathcal{I}_\mu^s(\mathcal{K}_\mu))^{\bot}$ satisfying $$\|\mathcal{R}(x_0, N, s)\|_{W^{-s}(\mathcal{K}_\mu)} \leq \frac{C_{s}}{N},$$ and for all $\mu \in spec(\Box)$ and $\mathcal{D}\in \mathcal{I}^s(\mathcal{B}_\mu^s)$, $$|c_\mathcal{D}| \leq C_{s}.$$
\end{theorem}
$\textit{Remark}:$  For sufficiently smooth functions, the sequence of values $\{\mathcal{D}(f)\}$ converges fast enough for the series $(\ref{155})$ to converge (see Section 7.3). 

\subsection{Acknowledgements}

This was written under Giovanni Forni for my PhD thesis, and I am sincerely grateful for his help.

\section{Orthogonal Bases for Models}

Let $A = \left(\begin{array}{rr}
a & b\\
c & d
\end{array}\right) \in SL(2, \mathbb{R}).$\\
For Casimir parameter $\mu > 0$, let $\nu = \sqrt{1 - \mu}$ be a representation parameter.  We denote by $\mathcal{H}_\mu$ the following models for the principal and complementary series representation spaces.  In the first model (the line model) the Hilbert space is a space of functions on $\R$ with the following norms.  If $\mu \geq 1$, then $\nu \in i \mathbb{R}$ and $\| f \|_{\mathcal{H}_\mu} = \| f \|_{L^2(\mathbb{R})}.$  If $0 < \mu < 1$, then $0 < \nu < 1$ and $$\| f \|_{\mathcal{H}_\mu} = \left(\int_{\mathbb{R}^2} \frac{f(x) \overline{f(y)}}{|x - y|^{1 - \nu}} dx dy\right)^{1/2}.$$  The group action is defined by $$\pi_\nu : SL(2, \mathbb{R}) \rightarrow \mathcal{U}(\mathcal{H}_\mu)$$ $$\pi_{\nu} (A) f(x) = |-c x + a|^{-(\nu + 1)} f(\frac{d x - b}{-c x + a}),$$ where $x\in \mathbb{R}$.  

By the change of variable $x = \tan(\theta)$, we have the circle models $\mathcal{H}_\mu = L^2([\frac{-\pi}{2}, \frac{\pi}{2}], \frac{d\theta}{\cos^2(\theta)})$ for the principal series, 
and $$\| f \|_{\mathcal{H}_\mu} = \left(\int_{[-\pi/2, \pi/2]^2} \frac{f(\tan\theta) f(\tan \theta')}{|\tan\theta - \tan \theta'|^{1 - \nu}} \frac{d\theta d\theta'}{\cos^2(\theta) \cos^2(\theta')}\right)^{1/2}$$ for the complementary series.

Computing derived representations, we get 
\begin{claim}\label{136}
Let $\mu > 0$. The vector fields for the $\mathcal{H}_\mu$ model on $\mathbb{R}$ are $$\begin{array}{lll}
X = -(1 + \nu) - 2 x \frac{\partial}{\partial x}; & \Theta = -(1 + \nu) x - (1 + x^{2}) \frac{\partial}{\partial x};\\
Y = (1 + \nu) x - (1 - x^{2}) \frac{\partial}{\partial x}; & U = -\frac{\partial}{\partial x}; \\
V = (1 + \nu) x + x^{2}\frac{\partial}{\partial x}.
\end{array}$$
By the change of variable $x = \tan(\theta)$, the vector fields in the circle model are:
$$\begin{array}{lll}
X = -(1 + \nu) - \sin (2\theta) \frac{\partial}{\partial \theta}; & \Theta = -(1 + \nu) \tan(\theta) - \frac{\partial}{\partial \theta}; \\
Y = (1 + \nu) \tan(\theta) - \cos(2 \theta) \frac{\partial}{\partial \theta}; & U = -\cos^{2} (\theta)\frac{\partial}{\partial \theta}; \\
V = (1 + \nu) \tan(\theta) + \sin^{2} (\theta)\frac{\partial}{\partial \theta}.
\end{array}$$
\end{claim}

 For $\mu \leq 0$, let $L_{hol}^2(H, d\lambda_\nu)$ be the upper half-plane model for the holomorphic discrete series, where $d\lambda_\nu := y^{\nu - 1} dx dy$ and $\nu = \sqrt{1 - \mu} \in \{2j - 1\}_{n \in \mathbb{Z}^+}$ is the representation parameter.  This model has the group action $\pi_\nu: SL(2, \mathbb{R})\rightarrow \mathcal{U}\left(L_{hol}^2(H, d\lambda_\nu)\right)$ defined by \begin{equation}\label{551} \pi_{\nu} (A): f(z) \rightarrow (-c z + a)^{-(\nu + 1)} f(\frac{d z - b}{-c z + a}).\end{equation}  The anti-holomorphic discrete series case occurs when $\nu = -\sqrt{1 - \mu} < 0$, but we only consider the holomorphic case because there is a complex anti-linear isomorphism between two series of the same Casimir parameter.  The space $L_{hol}^2(H, d\lambda_\nu)$ is said to be of $\textit{lowest weight } n := \frac{1 + \nu}{2}$.

The map $\alpha: D \rightarrow H : \xi \rightarrow - i\frac{\xi + 1}{\xi - 1}$ is a conformal map between $D$ and $H$.  For each $\nu \geq 1$, the unit disk model for the holomorphic discrete series is denoted $L_{hol}^2(D, d\sigma_\nu)$ and has the measure $d\sigma_\nu :=  4 \frac{(1 - |\xi|^2)^{\nu - 1}}{|\xi - 1|^{2(\nu + 1)}} du dv,$ which is calculated by change of variable. 

\begin{claim}\label{634}
Let $\mu \leq 0$.  Then the vector fields in $L_{hol}^2(H, d\lambda_\nu)$ are:
$$\begin{array}{lll}
X = -(1 + \nu) - 2 z \frac{\partial}{\partial z}; &
\Theta = -(1 + \nu) z - (1 + z^{2}) \frac{\partial}{\partial z} \\
Y = (1 + \nu) z - (1 - z^{2}) \frac{\partial}{\partial z}; &
U = -\frac{\partial}{\partial z}; \\
V = (1 + \nu) z + z^{2}\frac{\partial}{\partial z}.
\end{array}$$
By changing variables via the Mobius transformation $\alpha$, the vector fields for $L_{hol}^2(D, d\sigma_\nu)$ are:
$$\begin{array}{lll}
X =  -(1 + \nu) + (\xi^2 - 1)\frac{d}{d\xi}; & 
\Theta = (1 + \nu) i(\frac{\xi + 1}{\xi - 1}) - 2i\xi\frac{d}{d\xi}; \\
Y = -(1 + \nu) i(\frac{\xi + 1}{\xi - 1}) + i(\xi^2 + 1) \frac{d}{d\xi}; &
U = i\frac{(\xi - 1)^2}{2} \frac{d}{d\xi} \\
V = -(1 + \nu) i(\frac{\xi + 1}{\xi - 1}) + i \frac{(\xi +1)^2}{2} \frac{d}{d\xi}.
\end{array}$$
\end{claim}

We generate the basis for the principal and complementary series from a single element $u_0$ using the creation and annihilation operators $\eta_{\pm} = X \pm i Y$.   The basis for the discrete series is generated by applying the $\eta_+$ to the vector of lowest weight $n$.  

We calculate concrete formulas for the orthogonal basis vectors $\{u_k\}$ in Appendix A, and we present them here.  

\begin{lemma}\label{554}

(i) Let $\mu > 0$.  Then the set $\{u_k = e^{2i k \theta} \cos^{1 + \nu}(\theta)\}_{k \in \mathbb{Z}}$ is an orthogonal basis for the circle model of $\mathcal{H}_\mu$.  Moreover, if $\mu \geq 1$, then for all $k$, $$\| u_k\|_{\mathcal{H}_\mu}^2  = \pi.$$  If the Laplace-Beltrami operator $\triangle_M$ has a spectral gap $\mu_0 > 0$, then there is a constant $C_{SM} > 0$ such that for any $0 < \mu < 1$, $$C_{SM}^{-1} (1 + |k|)^{-\nu} \leq \| u_k \|_{\mathcal{H}_\mu}^2 \leq C_{SM} (1 + |k|)^{-\nu}.$$

(ii) Let $\mu \leq 0$ and $n = \frac{1 + \nu}{2}$ be the lowest weight.  Then $\{u_k = \left(\frac{z - i}{z + i}\right)^{k - n} \frac{1}{(z + i)^{\nu + 1}}\}_{k = n}^{\infty}$ is an orthogonal basis, and for all $k \geq n$, $$\| u_k\|_{L^2(H, d\lambda_\nu)}^2 = \frac{\pi}{\nu} 4^{-\nu} \left(\frac{(k - n)! \nu!}{(k + n - 1)!}\right).$$
\end{lemma}

By Claim $\ref{133}$, $\{u_k\}_{k}$ is an orthogonal basis for $\mathcal{H}_\mu$, $L^2(H, d\lambda_\nu)$ respectively.  The calculations concerning their norms are given in Appendix $A$.  We remark that in the complementary series, the values $\|u_k\|_{\mathcal{H}_\mu}$ converge to zero as $\mu \rightarrow 0^+$ (see also Lemma 2.1 of \cite{2}).  We assume throughout that $\triangle_{M}$ has a spectral gap in order that the values $\|u_k\|_{\mathcal{H}_\mu}$ have uniform lower bound in $\mu \in (0, 1).$  With this, we can prove a uniform constant in our estimates $(\ref{137})$ for the complementary series.

Sobolev norms of the basis vectors $\{u_k\}$ are given by \begin{equation}\| u_k\|_s^2 = (1 + \mu + 8k^2)^s \|u_k\|^2.\end{equation} (See Claim $\ref{133}$.)

\section{Relevant distributions in models}

\subsection{Invariant distributions}

Let $\mu \in spec(\triangle_M) - \{0\}$ and write $f = \Phi \cdot \cos^{1 + \nu} \in \mathcal{H}_\mu$ in circle coordinates, where $$\Phi(\theta) = \sum_{k = -\infty}^{\infty} c_k e^{2 i k \theta}.$$  Then define $$\delta^{(0)}(f) := \Phi(\frac{\pi}{2}),$$ and now formula (40) of \cite{2} shows $\delta^{(0)}$ is $U$-invariant, and Theorem 1.1 of \cite{2} proves \begin{equation}\label{2005}\delta^{(0)} \in W^{-((1 + \Re\nu)/2 + \epsilon)}(\mathcal{H}_\mu),\end{equation} for all $\epsilon > 0$

An important property of $\delta^{(0)}$ is that functions in $Ker(\{\delta^{(0)}\})$ decay at infinity.

\begin{lemma}\label{3011}
Let $\mu \in spec(\triangle_M) - \{0\}$, $\epsilon > 0$ and suppose $\triangle_M$ has a spectral gap $\mu_0 > 0$.  Then there is a constant $C_{s, \epsilon, SM} > 0$ such that for all $f\in W^{(1 + \Re\nu)/2 + 2 \epsilon}(\mathcal{H}_\mu)\cap Ann(\{\delta^{(0)}\})$ and $x\in\mathbb{R}$, we have $$|f(x)| \leq C_{\epsilon, SM} (1 + |x|)^{-(1 + \Re\nu + \epsilon)} \| f \|_{W^{(1 + \Re\nu)/2 + 2 \epsilon}(\mathcal{H}_\mu)}.$$
\end{lemma}

For $k\in \mathbb{Z}$, formally define the linear functional $\hat{\delta}_{k/T}$ by $$\hat{\delta}_{k/T}(f) = \int_\mathbb{R} f(x) e^{-2\pi i k/T x} dx = \hat{f}(\frac{k}{T})$$  When $0 < \mu < 1$, Sobolev embedding shows that sufficiently regular functions in $\mathcal{H}_\mu$ are also in $L^1(\R)$ (See Lemma $\ref{901}$).  

On the other hand, when $\mu \geq 1$, $C^\infty(\mathcal{H}_\mu)$ is not contained in $L^1(\R)$, nor is $L^1(\R)\cap C^{\infty}(\mathcal{H}_\mu)$ dense in $C^{\infty}(\mathcal{H}_\mu)$.  
Lemma\ $\ref{3011}$ gives that all functions that are sufficiently smooth and in $Ker(\delta^{(0)})$ are also in $L^1(\mathbb{R})$.  So for $\mu \geq 1$, we extend the definition of the Fourier transform $\mathcal{F}$ to any $f\in W^{(1 + \Re\nu)/2 + \epsilon}(\mathcal{H}_\mu)$ by setting \begin{equation}\label{224}\hat{f} := \mathcal{F}(f) = \mathcal{F}(f - \delta^{(0)}(f)\cos^{1 + \nu}\circ \arctan).\end{equation} 

\begin{lemma}\label{901}
Let $\mu \in spec(\triangle_M - \{0\}), \epsilon > 0, T > 0$ and $k\in \mathbb{Z}.$  Then $\hat{\delta}_{k/T} \in W^{-((1 + \Re\nu)/2 + \epsilon)}(\mathcal{H}_\mu)$ is a $\phi_T^U$-invariant distribution.
\end{lemma}
$\textit{Proof}:$
When $0 < \mu < 1$ and $f\in C^{\infty}(\mathcal{H}_\mu)$, write $f = \Phi\circ\arctan(x) \cos^{1 + \nu}\circ \arctan(x).$  Sobolev embedding gives $$|\hat{\delta}_{k/T}(f)| \leq \| \Phi\circ\arctan \|_{L^\infty(\R)} \int_\R  \cos^{1 + \nu}\circ \arctan(x) dx \leq C_{\epsilon, \nu} \| f \|_{W^{(1 + \Re\nu)/2 + \epsilon}(\mathcal{H}_\mu)}.$$  When $\mu \geq 1$, let $\tilde{f} := f - \delta^{(0)}(f)\cos^{1 + \nu}(\arctan x)$, and then the decay on $\tilde{f}$ given by Proposition $\ref{75}$ shows $$|\hat{\delta}_{k/T}(f)| \leq C_{\epsilon} \|f\|_{(1 + \Re\nu)/2 + \epsilon}.$$ 

 Additionally, the change of variable $x - T \mapsto x$ shows that $\hat{\delta}_{k/T}$ is $\phi_T^U$-invariant, $$\hat{\delta}_{k/T}f(\cdot - T) = \int_{\mathbb{R}} \tilde{f}(x - T) e^{-2\pi i x k/T}dx = \hat{\delta}_{k/T}(f). \ \ \Box$$  

The discrete series case is similar to the principal and complementary series cases.  Let $\mu \leq 0$ and $n = \frac{\nu + 1}{2}$ be the lowest weight.  By the change of variable $\xi = \left(\frac{z - i}{z + i}\right)$, the basis for $L_{hol}^2(H, d\lambda_\nu)$ written in the unit disk model $L_{hol}^2(D, d\sigma_\nu)$ is $\{u_k(\xi) = \xi^{k - n}(\xi - 1)^{\nu + 1}\}_{k = n}^{\infty}.$  Then any $f\in \mathcal{H}_\mu$ has the form $f = \Phi \cdot u_n \in \mathcal{H}_\mu$, where $$\Phi(\xi) = \sum_{k = n}^{\infty} c_k \xi^{k - n}.$$  Now define $\delta^{(0)}(f) := \Phi(1)$, so $\delta^{(0)} \in W^{-((1 + \nu)/2 + \epsilon)}(\mathcal{H}_\mu)$, again by formula (40) and Theorem 1.1 of \cite{2}.   

For $k \in \mathbb{Z}$, there are also distributions given by Fourier transforms of delta distributions along $\mathbb{R}\times \{i y\}.$  For $f\in W^{1 + \epsilon}(\mathcal{H}_\mu)$,  $k\in \mathbb{Z}$ and $y\in \mathbb{R}^+$, define $$\hat{\delta}_{k, y}(f) = \int_{\mathbb{R}} f(x + i y) e^{- 2\pi i k (x + i y)} dx.$$   

 \begin{lemma}\label{627}
Let $\mu \leq 0$, $k\in \mathbb{Z}$, $T > 0$ and $\epsilon > 0$.  Then $$\hat{\delta}_{k/T, y}\in W^{-(1 + \epsilon)}(H, d\lambda_\nu)$$ is a $\phi_T^U$-invariant distribution. 
\end{lemma}

Lemma $\ref{627}$ will follow immediately from Lemma $\ref{447}$, which proves functions in $W^{1 + \epsilon}(H, d\lambda_\nu)$ have some decay at infinity.  Moreover,

\begin{lemma}\label{767}
Let $\mu \leq 0$, $k \in \mathbb{Z}$ and $y_1, y_2 > 0$.  Then $\hat{\delta}_{k/T, y_1} = \hat{\delta}_{k/T, y_2},$ and if $k \leq 0$, then $\hat{\delta}_{k/T, y_1} = 0$.
\end{lemma}

Lemma $\ref{767}$ follows from Lemma $\ref{447}$ and Cauchy's theorem, and its proof is given in subsection $\ref{3020}$.  We therefore drop the subscript $y$ and write $$\hat{\delta}_{k/T} := \hat{\delta}_{k/T, y},$$ for $\textit{ any }$ $y > 0$. 

\subsection{Additional distributions at infinity}

The following distributions are$\textit{ not }\phi_T^U$-invariant and are introduced as a technical tool.\\

Let $\mu \in spec(\triangle_M - \{0\})$, and for all $r \in \mathbb{N}$, define $$\delta^{(r)} := (\Theta^r\delta^{(0)}).$$  Then Lemma 6.3 of Nelson, Analytic Vectors (\cite{8}) together with $(\ref{2005})$ shows  $$|\delta^{(r)}(f)| = |\delta^{(0)}(\Theta^r f)| \leq C_\nu \|\Theta^r f\|_{(1 + \Re\nu)/2 + \epsilon} \leq C_{\nu, r} \| f \|_{r + (1 + \Re\nu)/2 + \epsilon}.$$  Hence, $\delta^{(r)} \in W^{-(r + (1 + \Re\nu)/ 2 + \epsilon)}(\mathcal{H}_\mu).$ 

The$\textit{ Key Point }$ in proving our estimate for the cohomological equation is that functions that annihilate distributions at infinity have additional decay, and their derivatives do too.  We have the following stronger form of Lemma $\ref{3011}$.
  
\begin{proposition}\label{75}
Let $\mu \in spec(\triangle_M) - \{0\}$, $\epsilon > 0$ and $s \geq 0$, and suppose $\triangle_M$ has a spectral gap $\mu_0 > 0$.  Then there is a constant $C_{s, \epsilon, SM} > 0$ such that for all $f\in W^{s + (1 + \Re\nu)/2 + \epsilon}(\mathcal{H}_\mu)\cap Ann(\{\delta^{(r)}\}_{r = 0}^{\lfloor{s - 1}\rfloor})$, $x\in\mathbb{R}$ and integers $0 \leq r \leq s$, we have $$|f^{(r)}(x)| \leq C_{s, \epsilon, SM} (1 + |\nu|)^{r} (1 + |x|)^{-(s  + r + 1 + \Re\nu)} \| f \|_{W^{s + (1 + \Re\nu)/2 + \epsilon}(\mathcal{H}_\mu)}.$$
\end{proposition}

$\textit{Proof}:$ Write $f(\theta) = \Phi(\theta)\cos^{1 + \nu}(\theta)$.   Because $f$ is smooth in the representation theory sense, $\Phi$ has a Taylor series about $\frac{\pi}{2}$.  Moreover, because $f \in Ann(\{\delta^{(r)}\}_{r = 0}^{s - 1}$), we have $\Phi^{(r)}(\frac{\pi}{2}) = 0$ for all $0 \leq r \leq s - 1$.  So $\Phi$ decays, which forces $f$ to decay as well.  The estimate is a straightforward calculation and is similar to the proof of Proposition $\ref{929}$. $\ \ \Box$\\

Let $\mu \leq 0$.  Similarly, \begin{equation}\label{480}\delta^{(r)} := (\Theta^r \delta^{(0)}) \in W^{-(r + (1 + \nu)/2 + \epsilon)}(H, d\lambda_\nu).\end{equation} 

Recall that the parameter $\nu$ tends to infinity.  For fixed regularity $s$ and $\nu < s$, we estimate the transfer function in the same way that we did for the principal and complementary series.  We use the following proposition proven in Subsection $\ref{3020}$: 
\begin{proposition}\label{929}
Let $\mu \leq 0$, and let $r, s\in \mathbb{N}_0$ satisfy $0 \leq r < (s - 1)/2$ and $s \geq 4$.   Also let $f\in W^s(H, d\lambda_\nu) \cap Ann(\{\delta^{(r)}\}_{r = 0}^{\lfloor \frac{s - 1}{2}\rfloor - 1})$.  Then there is a constant $C_s > 0$ such that for all $z \in H$, $$|f^{(r)}(z)| \leq C_s \| f \|_s (1 + |z|)^{-(s/2 + \nu + r)}.$$
\end{proposition}
The case $s \geq \nu$ is different, and for this we do not use the additional distributions.

\section{Cohomological equation for the principal and complementary series}
Define $\mathcal{I}^s(\mathcal{H}_\mu) := \{\mathcal{D} \in W^{-s}(\mathcal{H}_\mu) | L_{-T} \mathcal{D} = 0\}$ and $Ann(\mathcal{I}^s(\mathcal{H}_\mu)) := \{f \in W^s(\mathcal{H}_\mu) | \mathcal{D} f = 0 \text{ for all } \mathcal{D}\in \mathcal{I}^s(\mathcal{H}_\mu)\}.$\\

$\textit{Throughout this section}$, let $\mu \in spec(\triangle_{M}) - \{0\}, T > 0$,  $r \geq 0$ and suppose $\triangle_{M}$ has a spectral gap $\mu_0 > 0$.  

The main theorem of this section is the following:
\begin{theorem}\label{165}
For all $f\in Ann(\mathcal{I}^{2r + 3/2}(\mathcal{H}_\mu))$, there exists a unique $\mathcal{H}_\mu$ solution $u$ to the cohomological equation $u \circ \phi_T^U - u = f.$
Additionally, there is a constant $C_{r, SM} > 0$ such that 
\begin{equation}\label{623}\| u \|_{W^r(\mathcal{H}_\mu)} \leq C_{r, SM, T} \| f \|_{W^{2r + 3/2}(\mathcal{H}_\mu)}.\end{equation}
\end{theorem}

$\textit{Remark}:$  We actually prove the tame estimate $$\| u \|_{W^r(\mathcal{H}_\mu)} \leq C_{r, SM} (1 + |\nu|)^r \| f \|_{W^{r + 3/2}(\mathcal{H}_\mu)}$$ in each irreducible component.  Because there exists infinitely many irreducible components, the representation parameters $\nu$ may tend to infinity, so we absorb $\nu$ using the Casimir operator $\Box$ and obtain $(\ref{623})$.

\subsection{Proof of Theorem $\ref{206}$}

The following theorem essentially proves Theorem $\ref{165}$.

\begin{theorem}\label{206}
If $f\in W^{2r + 3/2}(\mathcal{H}_\mu) \cap Ann(\{\hat{\delta}_{k/T}\}_{k = -\infty}^{\infty}\cup\{\delta^{(k)}\}_{k = 0}^{r + 1})$, then there exists a unique $\mathcal{H}_\mu$ solution $u$ to the cohomological equation \begin{equation}\label{92}u(x - T) - u(x) = f(x),\end{equation}
and there is a constant $C_{r, SM} > 0$ such that \begin{equation}\label{205}\| u \|_{W^r(\mathcal{H}_\mu)} \leq \frac{C_{r, SM}}{T} \ \| f \|_{W^{2r + 3/2}(\mathcal{H}_\mu)}.\end{equation}
\end{theorem}

To ease notation, define $$s(\nu, \epsilon): = s + (1 + \Re\nu)/2 + \epsilon.$$  

\begin{lemma}\label{49}
Let $s \geq 0$, and $f\in W^{s(\nu, \epsilon)}(\mathcal{H}_\mu) \cap Ann(\{\delta^{(r)}\}_{r = 0}^{s - 1}\cup \{\hat{\delta}_{k/T}\}_{k = -\infty}^{\infty}).$  Let $u$ be defined by $$u(x) = \sum_{n = 1}^{\infty} f(x + kT).$$  Then $u$ is a solution to $(\ref{92})$ and there is a constant $C_{s, \epsilon, SM} > 0$ such that for all $x\in \mathbb{R}$ and $0 \leq r \leq s$,  $$|u^{(r)}(x)| \leq \frac{C_{s, \epsilon, SM}}{T} (1 + |\nu|)^r (1 + |x|)^{-(s + r + \Re\nu)} \| f \|_{W^{s(\nu, \epsilon)}(\mathcal{H}_\mu)}.$$
\end{lemma}
$\textit{Proof}:$
Then $u$ is formally a solution.  When $x \geq 0$, Proposition $\ref{75}$ gives $$|u^{(r)} (x) \cdot T| = |\frac{d^r}{dx^r} \sum_{n = 1}^{\infty} f(x + n T) \cdot T| \leq \sum_{n = 1}^{\infty} |f^{(r)}(x + n T)| \cdot T$$ $$ \leq C_{s, \epsilon, SM} (1 + |\nu|)^r \| f \|_{s(\nu, \epsilon)} \sum_{n = 1}^{\infty} (|x + n T| + 1)^{-(s + r + 1 + \Re\nu)} \cdot T.$$  Then by the integral estimate we conclude.  

For $x < 0$, Proposition $\ref{75}$ shows $f \in L^1(\mathbb{R})$, so the Poisson summation formula gives $$\sum_{n\in\mathbb{Z}} f(x + n T) = \frac{1}{T}\sum_{n\in \mathbb{Z}} \hat{f}(\frac{n}{T}) e^{2\pi i x n/ T} = 0.$$ Therefore, $u(x) = \sum_{n = 1}^{\infty} f(x + n T) = -\sum_{n = 0}^{\infty} f(x - n T).$  Again, Proposition $\ref{75}$ and the integral estimate proves the lemma. $\ \ \Box$\\
 
 $\textit{Proof of Theorem } \ref{206}:$
 From the definition of $\triangle$ and the principal and complementary series formulas for the vector fields $X, Y$ and $\Theta$ given in Claim $\ref{136}$, one finds $\triangle^r$ consists of terms of the form $(1 + \nu)^m x^k \frac{d^j}{dx^j}$, where $0 \leq k - j \leq r$ and $j + m = r$.
 
Then $$\| u\|_{W^r(\mathcal{H}_\mu)} \leq \sum_{0 \leq j, k - j \leq r} (1 + |\nu|)^{r - j} \|(|x| + 1)^k u^{(j)}\|_{\mathcal{H}_\mu}$$ \begin{equation}\label{976}\leq  \frac{C_{s, \epsilon, SM}}{T} (1 + |\nu|)^r \| f\|_{s(\nu, \epsilon)} \|(|x| + 1)^{r - (s + \Re\nu)} \|_{\mathcal{H}_\mu}.\end{equation}
 
 Note that if $\mu \geq 1$, then $\| \cdot\|_{\mathcal{H}_\mu} = \|\cdot\|_{L^2(\mathbb{R})}$ and $\Re\nu = 0$.  So for all $0 \leq r < s - 1/2$, $(\ref{976}) < \infty.$  In this case, $s(\nu, \epsilon) = s + \frac{1}{2} + \epsilon$, so $$\|u\|_{r} \leq  \frac{C_{s, \epsilon, SM}}{T} (1 + |\nu|)^r \| f\|_{s + 1/2 + \epsilon}.$$  In particular, this holds for $r = s - 1/2 - \epsilon$.  Then replacing $s$ with $r + 1/2 + \epsilon$ proves Theorem $\ref{206}$ for $\mu \geq 1$ and $r\in \mathbb{N}_0$.  
 
 When $0 < \mu < 1$, then $\Re \nu = \nu$ in formula $(\ref{976})$.  Lemma $\ref{1095}$ shows that \begin{equation}\label{2050}\| g \|_{\mathcal{H}_\mu} \leq C_\nu \left(\|g\|_{L^\infty(\mathbb{R})} + \|g\|_{L^1(\mathbb{R})}\right).\end{equation}  There are only finitely many values in $spec(\Box) \subset (0, 1)$, so the constant $C_\nu$ in $(\ref{2050})$ satisfies $C_\nu \leq C$ for some absolute constant $C.$ 
 
 Observe that for all $0 \leq r < s + \nu - 1$, $(|x| + 1)^{r - (s + \nu)}\in L^1(\mathbb{R}) \cap L^{\infty}(\mathbb{R}).$  As before, this holds for $r = s + \nu  - 1 - \epsilon$.  Setting $s = r + 1 + \epsilon$, we see $$s(\nu, \epsilon) = r + 1 - \nu + (1 + \nu)/2 + 2\epsilon \leq r + 3/2.$$  This proves the estimate in Theorem $\ref{206}$ for $r\in \mathbb{N}_0$.  

Let $0 \leq s \leq r$, so that we have the continuous injections $W^r(\mathcal{H}_\mu) \subseteq W^s(\mathcal{H}_\mu) \subseteq \mathcal{H}_\mu.$  The operator $(I + \triangle)$ is an essentially self-adjoint, strictly positive operator on $\mathcal{H}_\mu$, so for any $\alpha\in [0, 1]$, the operator $(I + \triangle)^{r \alpha + s (1 - \alpha)}$ on $\mathcal{H}_\mu$ is defined by the spectral theorem.  Then define $W^{r \alpha + s (1 - \alpha)}$ $(\mathcal{H}_\mu) := dom((I + \triangle)^{r \alpha + s (1 - \alpha)})_{|_{\mathcal{H}_\mu}}$ and give it the norm $$\| f \|_{r \alpha + s (1 - \alpha)} := \| (1 + \triangle)^{r \alpha + s (1 - \alpha)} f \|_{\mathcal{H}_\mu}.$$  In this sense, $\{W^{r}(\mathcal{H}_\mu)\}_{r\in \R_{\geq 0}}$ is an interpolation family.  Then because the estimate in Theorem $\ref{206}$ holds for all integers $r \geq 0$, the interpolation theorem given at Theorem 5.1 in \cite{11} completes the proof.  $\ \ \Box$  
\\ 

\subsection{Proof of Theorem $\ref{165}$}

\begin{theorem}\label{207}
Theorem $\ref{206}$ holds under the weakened hypothesis that $$f\in W^{r + 3/2}(\mathcal{H}_\mu) \cap Ann\left(\{\hat{\delta}_{n/T}\}_{n = -\infty}^{\infty}\cup\{\delta^{(0)}\}\right).$$
\end{theorem}

To begin, set $\chi_0 := u_0 ( = \cos^{1 + \nu}(\cdot))$, and recursively define $\{\chi_k\}_{k = 1}^r$ by  \begin{equation}\label{161}\chi_{k + 1} := (\chi_{k}\circ \phi_T^{U} - \chi_k).\end{equation}   Then $\chi_k$ is a coboundary for all $k \geq 1$.  We show $\{\chi_k\}_{k = 1}^r$ is a basis in the dual space to $\langle\{\delta^{(k)}\}_{k = 0}^{r}\rangle$ and obtain a bound for each $\|\chi_k\|_{W^r(\mathcal{H}_\mu)}.$  For this, we study the distributions $\phi_T^U\delta^{(k)}$.

A calculation based on Claim $\ref{133}$, parts $ii)$ and $iii)$, proves
\begin{lemma}\label{560}
Let $\mu\in spec(\Box)$, $r\geq 0$ and $f\in C^{\infty}(\mathcal{H}_\mu)$.
  If $r$ is even then $$\mathcal{L}_U \delta^{(r)} = \frac{1}{2}\sum_{j = 0}^{[\frac{r}{2}] - 1}\left(2 (\nu + 1)(2i)^{2j}(_{2 j + 1}^{r}) - (2i)^{2(j + 1)}(_{2 (j + 1)}^{r})\right)\delta^{(r - 2j - 1)},$$ and if $r$ is odd, then 
$$\mathcal{L}_{U} \delta^{(r)} = \frac{i}{2}(\nu + 1)(2i)^r\delta^{(0)}$$ $$ + \frac{1}{2}\sum_{j = 0}^{[\frac{r}{2}] - 1} \left(2 (\nu + 1)(2i)^{2j}(_{2 j + 1}^{r}) - (2i)^{2(j + 1)}(_{2 (j + 1)}^{r})\right)\delta^{(r - 2j - 1)}. \ \ \Box$$
\end{lemma}

Lemma $\ref{560}$ gives coefficients $\{c_{j, k}\}_{0 \leq j, k \leq r}\subset \mathbb{C}$ such that ${\mathcal{L}_U}|_{\langle\{\delta^{(k)}\}_{k = 0}^r\rangle} = (c_{j, k})_{j, k}$ is an $r\times r$ strictly upper triangular matrix.  

Exponentiating, we get coefficients $\{e_{j, k}\}_{0 \leq j, k \leq r}\subset \mathbb{C}$ such that \begin{equation}\label{2052}\phi_{-T}^U|_{\langle\{\delta^{(k)}\}_{k = 0}^r\rangle} = \left(e_{j, k}\right)_{j, k}\end{equation} is an $r\times r$ upper triangular matrix where  $e_{j, k} = \frac{(-T)^{k - j}}{(k - j)!} a_{j, k}$ for some coefficients $\{a_{j, k}\} \subset \mathbb{C}$ and $e_{j j} = 1$ for all $j$.  
\begin{lemma}\label{1100}
Let $\{\chi_k\}_{k = 0}^{r}$ be defined by $(\ref{161})$.  Then for all $1 \leq j\leq k \leq r$, \begin{equation}\label{260}\delta^{(j)}(\chi_k) = \Bigg\{\begin{array}{rr} \Pi_{j = 0}^{k - 1}e_{j, j + 1}\text{ if } j = k\\ 0 \ \ \ \text{ if } j < k\end{array}\end{equation}
\end{lemma}
$\textit{Proof}:$
This follows by induction using the identity $(\phi_{-T}^U\delta^{(j)})(\chi_{k}) = \sum_{m = 0}^{j} e_{m, j} \delta^{(m)}(\chi_{k})$ from $(\ref{2052}). \ \ \Box$

For convenience, we define $\Pi_k := \Pi_{j = 0}^{k - 1}|e_{j, j + 1}|$ for all $k \geq 1$.
\begin{lemma}\label{933}
Let $\mu\in spec(\Box)$, $\epsilon > 0$ and $f\in W^r(\mathcal{H}_\mu) \cap Ann(\{\delta^{(0)}\})$.  Then there are coefficients $\{\omega_{k, f}\}_{k = 1}^{r}$ and a constant $C_r > 0$ such that $$f_d :=  f - \sum_{k = 1}^r \omega_{k, f} \chi_k$$ is in $Ann(\{\delta^{(k)}\}_{k = 0}^r)$ and for all $1 \leq k \leq r$, \begin{equation}\label{978}|\omega_{k, f}| \leq \bigg\{\begin{array}{rr} \frac{C_r }{T^{k(k + 1)/2}} \| f \|_{k + (1 + \Re\nu)/2 + \epsilon} \text{ if } T < 1\\ C_r \| f \|_{k + (1 + \Re\nu)/2 + \epsilon} \text{ otherwise}.\end{array}\end{equation}
\end{lemma}
$\textit{Proof:}$
To ease notation, we will write $\omega_j := \omega_{j, f}$ for all $j$.  Recursively define $\omega_1 := \frac{\delta^{(1)}(f)}{\Pi_1}$ and if $\omega_j$ have been defined for $1 \leq j < k \leq r$, define $$\omega_k :=  \frac{\delta^{(k)}(f) - \sum_{j = 1}^{k - 1} \omega_j \delta^{(k)}(\chi_{j})}{\Pi_{k}}.$$  We first prove by induction that for all $0 \leq k \leq r$, $\delta^{(k)}(f_d) = 0.$  By assumption $\delta^{(0)}(f) = 0$, and because all functions in $\{\chi_k\}_{k = 1}^r$ are coboundaries, the flow invariance of $\delta^{(0)}$ implies $\delta^{(0)}(\chi_k) = 0$ for $k \geq 1$.  So $\delta^{(0)}(f_d) = 0.$  Now assume that $\delta^{(j)}(f_d) = 0$ for $0 \leq j < k$.  Moreover, by construction and Lemma $\ref{1100}$, $$\delta^{(k)}(f_d) = \delta^{(k)}(f) - \sum_{j = 1}^r \omega_j \delta^{(k)}(\chi_j)$$ $$ = \delta^{(k)}(f) - \sum_{j = 1}^k \omega_j \delta^{(k)}(\chi_j) =  \delta^{(k)}(f) - \sum_{j = 1}^{k - 1} \omega_j \delta^{(k)}(\chi_j) - \omega_k \Pi_k = 0,$$ from the definition of $\omega_k$.

For the estimate, we exponentiate the matrix $\mathcal{L}_U|_{\langle\{\delta^{(k)}\}_{k = 0}^r\rangle}$ in Lemma $\ref{560}$ and get $$e_{j, j + 1} = - T c_{j, j + 1} = - T (j + 1) [j + (\nu + 1)].$$ 
 It follows that for $\mu \in spec(\Box),$ there is a constant $C_r > 0$ such that $\Pi_{r} \geq C_r T^r.$

Now we prove by induction that $(\ref{978})$ holds for all $1 \leq k \leq r$.  Consider the case $T < 1$, and then the case $T \geq 1$ will be clear.  Recall from formula $(\ref{480})$ that $\delta^{(j)} \in W^{-(j + (1  + \nu)/2)}(\mathcal{H}_\mu)$.  Notice that $|\omega_1| \leq \frac{C}{T} \| f \|_{1 + (1 + \Re\nu)/2 + \epsilon}.$  Then assuming $(\ref{978})$ holds for $1 \leq j < k$, observe $$|\omega_k| \leq \frac{C_k}{T^k} \left( C_k \| f \|_{k + (1 + \Re\nu)/2 + \epsilon} +  \sum_{j = 1}^{k - 1} \frac{C_j \| f \|_{j + (1 + \Re\nu)/2}}{T^{j(j + 1)/2}} \delta^{(k)}(\chi_j)\right)$$ $$ \leq \frac{C_k}{T^{k(k + 1)/2}} \| f \|_{k + (1 + \Re\nu)/2 + \epsilon}. \ \ \Box$$ 

$\textit{Proof of Theorem } \ref{207}:$ Let $f_d$ and $\{\omega_k\}_{k = 1}^r$ be defined as in Lemma $\ref{933}$.  Because $f_d \in Ann(\{\delta^{(j)}\}_{j = 0}^{r})$, Theorem $\ref{206}$ shows that $f_d$ has a transfer function $u_d$ and there is a constant $C_r > 0$ such that $$\|u_d\|_{W^r(\mathcal{H}_\mu)} \leq C_r \| f_d \|_{2r + 3/2}.$$  

For each $1 \leq k \leq r$, $\chi_k$ is a coboundary by construction, and there is a constant $C_{r} > 0$ such that $\|\chi_{k}\|_{W^{r}(\mathcal{H}_\mu)} \leq C_{r}.$  Then define $$u := u_d + \sum_{k = 1}^r \omega_{k, f} \chi_{k - 1},$$ and one checks that $u(\cdot - T) - u(\cdot) = f.$  
Moreover, $$\|u\|_{r} \leq \|u_d\|_{r} + \sum_{k = 0}^{r - 1}|\omega_{k, f}| \|\chi_k\|_r$$ $$ \leq \frac{C_{r, SM}}{T}\left(\| f_d \|_{2r + 3/2} + \frac{1}{T^{r(r + 1)/2}} \| f \|_{r + 1}\right) \leq C_{r, T, SM} \| f \|_{2r + 3/2}.$$ 
The case when $T \geq 1$ follows in the same way.  Finally, $u$ is the unique $\mathcal{H}_\mu$ solution, because if $w$ is any $\mathcal{H}_\mu$ solution to $(\ref{92})$, then $w - u \in \mathcal{H}_\mu$ and is $T$-periodic, which means $w = u$ in $\mathcal{H}_\mu$. $\ \ \Box$\\

$\textit{Proof of Theorem }\ref{165}:$ This now follows from Theorem $\ref{207}$ by showing that $\mathcal{I}^{2r + 3/2}(\mathcal{H}_\mu)$ is precisely $S_0 := \langle\{\hat{\delta}_{n/T}\}_{n\in \mathbb{Z}}\cup \{\delta^{(0)}\}\rangle.$  Section 3 shows $S_0 \subset \mathcal{I}^{2r + 3/2}(\mathcal{H}_\mu)$.  For the other inclusion, suppose there exists $\mathcal{D}\in \mathcal{I}^{2r + 3/2}(\mathcal{H}_\mu) - S_0$.  Then let $f\in C^{\infty}(\mathcal{H}_\mu)\cap Ann(S_0)$ be such that $\mathcal{D}(f) \neq 0.$  By Theorem $\ref{207}$, $f$ has a smooth transfer function in the domain of $\mathcal{D}$, so $\mathcal{D}(f) = 0.$  Contradiction. $ \ \ \Box$ \\

\section{Cohomological equation for the discrete series}

The spaces $\mathcal{I}^{s}(H, d\lambda_\nu)$ and $Ann(\mathcal{I}^{s}(H, d\lambda_\nu))$ are defined analogously to $\mathcal{I}^s(\mathcal{H}_\mu)$ and $Ann(\mathcal{I}^s(\mathcal{H}_\mu))$.  The main theorem of this section is 
\begin{theorem}\label{357}
Let $\mu \leq 0$, $T > 0$, $r \geq 0$, and $f\in Ann(\mathcal{I}^{3r +4}(H, d\lambda_\nu)).$  Then there is a unique $L^2(H, d\lambda_\nu)$ transfer function $u$ satisfying $u\circ \phi_T^U - u = f,$ and there is a constant $C_{r, T} > 0$ such that $$\| u \|_{W^r(\mathcal{H}_\mu)} \leq C_{r, T}
\| f \|_{W^{3r + 4}(\mathcal{H}_\mu)}.$$
\end{theorem}

According our model, $$u\circ \phi_T^U - u = f \text{ means } u(\cdot - T) - u = f,$$ for $u, f \in L^2(H, d\lambda_\nu)$.  Throughout this section we will use the biholomorphic map $\alpha : D \rightarrow H : \xi \rightarrow -i\left(\frac{\xi + 1}{\xi - 1}\right) := z.$ 

We remind the reader that it suffices to only consider the holomorphic discrete series.  The argument is divided into two dissimilar cases, when $\nu < s$ and when $\nu \geq s$.  When $\nu < s$, the function $f$ does not have enough decay to easily estimate its transfer function, so we use the additional distributions at infinity as we did in our estimate for the principal and complementary series.  We do not use them when $\nu \geq s$.

\subsection{Proof of Theorem $\ref{546}$}

Our immediate goal is to prove

\begin{theorem}\label{546}
Let $\mu \leq 0$, $T > 0$, $r \geq 0$, and $f\in W^{3r + 4}(H, d\lambda_\nu)\cap Ann(\{\delta_{k/T}\}_{k\in \mathbb{Z}^+} \cup \{\delta^{(k)}\}_{k = 0}^{r + 1}).$  Then there is a constant $C_{r, T} > 0$ and a unique $L^2(H, d\lambda_\nu)$ transfer function such that for all $z\in H$, \begin{equation}\label{359}u(z - T) - u(z) = f(z),\end{equation} and $$\| u \|_{W^r(H, d\lambda_\nu)} \leq C_{r, T} 
\| f \|_{W^{3r + 4}(H, d\lambda_\nu)}.$$
\end{theorem}

Let $\tilde{s} := \lfloor\frac{s - 1}{2}\rfloor.$  Our method of proving this is the same as for the principal and complementary series. 

\begin{lemma}\label{529}
Let $\mu \leq 0$, and let $r, s$ be integers that satisfy $0 \leq r < \tilde{s}$ and $s \geq 4$.   Also let $T > 0$, $f\in W^s(H, d\lambda_\nu)\cap Ann(\{\hat{\delta}_{k/T}\}_{k\in \mathbb{Z}^+}\cup \{\delta^{(r)}\}_{r \geq 0}^{\tilde{s}})$.  Then there is a constant $C_s > 0$ and a unique $L^2(H, d\lambda_\nu)$ solution $u$ to the cohomological equation $u(z - T) - u(z) = f(z)$ such that for all $z\in H$, $$|u^{(r)}(z)| \leq \frac{C_{r, s}}{T} \nu^r \| f \|_{W^s(H, d\lambda_\nu)} (1 + |z|)^{-(s/2 + \nu + r - 1)}.$$ 
\end{lemma}
$\textit{Proof}:$  This time we use Proposition $\ref{929}$ and conclude as in Lemma $\ref{49}$. $\ \ \Box$.\\

$\textit{Proof of Theorem }\ref{546}:$ As in the proof of Theorem $\ref{206}$, one finds \begin{equation}\label{602}\| u\|_{r} \leq \frac{C_{r, s}}{T} \nu^r \| f \|_s \left(\int_{H} (1 + |z|)^{2r + 2 - s - 2 \nu} \Im(z)^{\nu - 1}dxdy\right)^{1/2}.\end{equation}  
Observe $(\ref{602}) < \infty$ whenever $2r + 1 - s - \nu < -2$, which holds whenever $2r + 2 < s$.  Then choose $s = 2r + 4$.  

Finally, Claim $\ref{133}$ shows that $\Box f = (1 - \nu^2) f$, so that $\nu^2 f = (1 - \Box) f$.  Then for $r \in 2\mathbb{N}_0$, $$\| \nu^r f \|_{2r + 4} = \|(1 - \Box)^{r/2} f\|_{2r + 4} \leq C_r \| f \|_{3r + 4},$$ by Lemma $6.3$ of Nelson \cite{8}.  Finally, interpolation gives the estimate for all $r \in \mathbb{R}^+$.  The solution is unique for the same reason as in Lemma $\ref{529}$.  $\ \ \Box$

Now we remove the additional distributions.
\begin{theorem}\label{218}
Theorem $\ref{546}$ holds under the weakened hypothesis that $f\in W^{3r + 4}(H, d\lambda_\nu)\cap Ann(\{\delta_{n/T}\}_{n\in \mathbb{Z}} \cap\{\delta^{(0)}\})$.
\end{theorem}

 Set $\chi_0 := u_n$ and given $\chi_k$, define $\chi_{k + 1} := \chi_k\circ \phi_{-T}^U - \chi_k.$  Lemmas $\ref{1100}$ and $\ref{933}$ do not depend on the particular representation space, and the same argument used in the proof of Theorem $\ref{207}$ now proves Theorem $\ref{218}$. $ \ \ \Box$
 
\subsection{Proof of Theorem $\ref{357}$}

$\textit{Throughout this subsection}$ let $\mu \leq 0, \nu \geq 5, r \geq 0, \nu \geq 3r + 4$ and $T > 0$.  Let $H^+ = \{z \in H | \Re z > 0\}$ and $H^- = \{z \in H | \Re z < 0\}$.  

For this case $\nu \geq 3r + 4$, we do not use any distributions at infinity. 

\begin{theorem}\label{736}
Let $f\in W^{3r + 4}(H, d\lambda_\nu) \cap Ann(\{\hat{\delta}_{k/T}\}_{k = 1}^{\infty}).$  Then there is a constant $C_r > 0$ and a unique $L^2(H, d\lambda_\nu)$ transfer function $u$ to the cohomological equation $(\ref{359})$, which satisfies $$\| u \|_{W^r(H, d\lambda_\nu)} \leq \frac{C_{r}}{T} \| f \|_{W^{3r + 4}(H, d\lambda_\nu)}.$$
\end{theorem}

Let $\{u_k\}_{k \geq n} \subset L^2(H, d\lambda_\nu)$ be the basis discussed in Section 2, and for $k \geq 0$, write $$u_{k}(z) := \left(\frac{z - i}{z + i}\right)^{k - n} \left(\frac{1}{z + i}\right)^{\nu + 1}.$$  

\begin{lemma}\label{726}
Let $s\in \mathbb{N}$ and satisfy $2r + 4 \leq s$.  Then $$u(z) := \sum_{m = 1}^{\infty} f(z + mT)$$ is a solution to the cohomological equation \begin{equation}\label{727} u(z - T) - u(z)  = f(z),\end{equation} and $u^{(r)}$ is in $L^2(H, d\lambda_\nu)$.
\end{lemma} 
$\textit{Proof}:$
Let $f(z) = \sum_{k = n}^{\infty} c_k u_k(z)$, and define $u(z) := \sum_{m = 1}^{\infty} f(z + mT).$  Then $u$ is formally solution to $(\ref{727})$.  Observe the basis elements $u_k$ decay like $(1 + |z|)^{-(1 + \nu)}$.  Because $\nu \geq 3r + 4$, Lemma $\ref{447}$ shows that if $z \in H$, then $$|f^{(r)}(z)| \leq C_{r, T, \nu, y} \| f \|_s \frac{1}{(1 + |z|)^{2}}. \ \ \Box$$

\begin{lemma}\label{679}
For all $k \geq 0$, $$u_{k + n}^{(r)}(z) = \sum_{j = 1}^{r}\sum_{l = 0}^j \tilde{c}_{j, r} \frac{k!}{(k - (j - l))!} \frac{(\nu + 1)!}{(\nu + 1 - j)!} u_{k + n + l - j}(z) (z + i)^{-r},$$ where we set \begin{equation}\label{693}\frac{k!}{(k - (j - m))!} := 0\end{equation} if $k < j - m$.
\end{lemma}
$\textit{Proof}:$
  Let $\alpha: D \rightarrow H : \xi \rightarrow -i\left(\frac{\xi + 1}{\xi - 1}\right)$ be the M$\ddot{\text{o}}$bius transformation from Section 2.  Switching to unit disk coordinates and then using formula $(\ref{350})$, there are constants $\{c_{j, r}\}_{j = 0}^r \subset \mathbb{C}$ such that \begin{equation}\label{678}u_{k + n}^{(r)}(z) = U^r(u_{k + n}\circ \alpha)(\xi) = \sum_{j = 1}^{r}c_{j, r}(\xi - 1)^{r + j} (u_{k + n}\circ\alpha)^{(j)}(\xi).\end{equation}    
A calculation shows that for all $k \geq 0,$ $u_{k + n} \circ \alpha(\xi) = \xi^k \left(\frac{1 - \xi}{-2i}\right)^{\nu + 1}.$  Using the notation in $(\ref{693})$, it follows that$$(u_{k + n}\circ \alpha)^{(j)}(\xi) = \sum_{l = 0}^j (_l^j) \frac{k!}{(k - (j - l))!}\xi^{k- (j - l)} \frac{(\nu + 1)!}{(\nu + 1 - l)!} (\xi - 1)^{\nu + 1 - l}.$$  Combining this with $(\ref{678})$, we conclude.  $ \ \ \Box$\\

Given $k, q, j \in \mathbb{N}_0$, define $$v_{k + n, q, j, T}(z) := \sum_{m = 1}^{\infty} (1 + |z + m T|)^q |u_{k + n}^{(j)}(z + m T)|.$$ 

\begin{lemma}\label{725}
Let $s \in \mathbb{N}$ be such that $2r + 4 \leq s$, and let $f\in W^s(\mathcal{H}_\mu)$.  Then there is a constant $C_s > 0$ such that $$\| u\|_{W^r(H^+, d\lambda_\nu)} \leq C_s \| f \|_s \sum_{0 \leq j, q - j \leq r} \nu^{r - j} \left(\sum_{k = 0}^{\infty} \| u_{k + n} \|_{W^s(H, d\lambda_\nu)}^{-2} \|v_{k + n, q, j, T}\|_{L^2(H^+, d\lambda_\nu)}^2\right)^{1/2}.$$
\end{lemma}
$\textit{Proof}:$
As in the proofs of Theorems $\ref{206}$ and $\ref{546}$, there exists $C_r > 0$ such that $$\| (1 + \triangle)^r u\|_{L^2(H^+, d\lambda_\nu)} \leq C_r \sum_{0 \leq j , q - j \leq r} \nu^{r - j} \|(1 + |z|)^q u^{(j)}(z)\|_{L^2(H^+, d\lambda_\nu)}.$$ For $z \in H^+$, we have $$|(1 + |z|)^q u^{(j)}(z)| \leq \sum_{k = 0}^{\infty} \sum_{m = 1}^{\infty} |c_{k + n}| (1 + |z + m T|)^q |u_{k + n}^{(j)}(z + m T)| $$ $$ = \sum_{k = 0}^{\infty} \left(|c_{k + n}| \| u_{k + n}\|_s\right) \left(\| u_{k + n}\|_s^{-1} \sum_{m = 1}^{\infty}  (1 + |z + mT|)^q |u_{k + n}^{(j)}(z + m T)|\right)$$ $$\leq \| f \|_s \left(\sum_{k = 0}^{\infty} \| u_{k + n}\|_s^{-2} \left(\sum_{m = 1}^{\infty} (1 + |z + mT|)^q |u_{k + n}^{(j)}(z + m T)|\right)^2\right)^{1/2}.$$  Now taking $L^2$ norms gives Lemma $\ref{725}$. $ \ \ \Box$ 
  
  An elementary calculation proves
  \begin{lemma}\label{743}
With assumptions as in Lemma $\ref{725}$, there exists a constant $C_j > 0$ such that $$\| v_{k + n, q, j, T}\|_{L^2(H^+, d\lambda_\nu)} \leq \frac{C_j}{T} \frac{(\nu + 1)!}{(\nu + 1 - j)!} \frac{k!}{(k - j)!} \frac{1}{\sqrt{\nu - s + 1} \cdot 2^{\nu - s}}\left(\frac{(k - j)! (\nu - s)!}{(k - j + \nu - s)!}\right)^{1/2}.$$
 \end{lemma}

\begin{proposition}\label{744}
 With assumptions as in Lemma $\ref{725}$, there is a constant $C_{s, j} > 0$ such that $$\sum_{k = 0}^{\infty} \| u_{k + n} \|_s^{-2} \|v_{k + n, q, j, T}\|_{L^2(H^+, d\lambda_\nu)}^2 \leq \frac{C_{j, s}}{T}.$$
\end{proposition}
$\textit{Proof}:$
Lemma $\ref{554}\ ii)$ together with Claim $\ref{133}\ iii)$ give \begin{equation}\label{4015}\| u_{k + n} \|_s^{-2} = \frac{\nu}{\pi} \cdot 4^{\nu} (1 + \mu + 8(k + n)^2)^{-s} \left(\frac{(k + \nu)!}{k! \nu!}\right).\end{equation}  Multiplying this by $\|v_{k + n, q, j, T}\|_{L^2(\mathcal{H}_\mu)}$, the factor $4^\nu$ is canceled by $2^{-2(\nu - s)}$, and growth from the factor $\left(\frac{k!}{(k - j)!} \right)^2 \left(\frac{(k + \nu)!}{k! \nu!}\right)$ is canceled by $(1 + \mu + 8(k + n)^2)^{-s}\left(\frac{(k - j)! (\nu - s)!}{(k - j + \nu - s)!}\right)$.  The proposition follows. $\ \ \Box$

$\textit{Proof of Theorem }\ref{736}:$  
Combining Proposition $\ref{744}$ with Lemma $\ref{725}$ proves $\|u\|_{W^r(H^+, d\lambda_\nu)}$ $\leq$ $\frac{C_{r, s}}{T} \nu^r \| f \|_s.$  To estimate $u$ on $H^-$, we let $y > 0$.  Because $s > 1$, Proposition $\ref{748}$ proves $f(\cdot  + i y) \in L^1(\mathbb{R})$, so the Poisson summation formula applies.  As in the proof of Lemma $\ref{49}$,
 $f\in Ann(\{\hat{\delta}_{k/T}\}_{k \geq 1}) = Ann(\{\hat{\delta}_{k/T}\}_{k = -\infty}^{\infty})$, and we conclude $$u(z) = \sum_{m = 1}^{\infty} f(z + m T) = -\sum_{m = 0}^{\infty} f(z - m T).$$  The same argument used to estimate $\| u \|_{W^r(H^+, d\lambda_\nu)}$ proves there is a constant $C_s > 0$ such that $\|u\|_{W^r(H^-, d\lambda_\nu)} \leq \frac{C_{r, s}}{T} \nu^r \| f \|_s.$  Then combining the estimates for $H^+$ and $H^-$, and setting $s = 2r + 4$ proves $$\| u\|_r \leq \frac{C_r}{T} \nu^r \| f \|_{2r + 4} \leq \frac{C_r}{T} \| f \|_{3r + 4}$$ when $r\in \mathbb{N}_0$.  The estimate for $r \geq 0$ and real follows by interpolation.    

Finally, the solution $u$ is unique for the same reason discussed at the end of Lemma $\ref{529}$, which completes the proof of Theorem $\ref{736}$. $\ \ \Box$
 
 $\textit{Proof of Theorem 5.1:}$ It remains to show the space of invariant distributions $\mathcal{I}^{3r + 4}(H, d\lambda_\nu)$ is modeled by $$S_0 := \Big\{\begin{array}{rr}\langle\{\hat{\delta}_{k/T}\}_{k \in \mathbb{Z}^+} \cup \{\delta^{(0)}\}\rangle \text{ if } \frac{1 + \nu}{2} < 3r + 4 \\ \langle\{\hat{\delta}_{k/T}\}_{k \in \mathbb{Z}^+}\}\rangle \text{ otherwise}.\end{array}$$  By Section 3, $S_0 \subset \mathcal{I}^{3r + 4}(\mathcal{H}_\mu),$ and the other inclusion follows from definitions and Theorems $\ref{218}$ and $\ref{736}$ (See, for example, the proof of Theorem $\ref{165}$).  $\ \ \Box$

\section{Proof of Theorem $\ref{302}$}

We prove Theorem $\ref{302}$, which states that $\tilde{\mathcal{I}}(SM) := \mathcal{I}(SM) - \int_{\oplus_\{\mu < 0\}} \langle \{\mathcal{D}_\mu^0\}\rangle d\beta(\mu)$ is the space of distributional obstructions to the existence of $L^2(SM)$ solutions for coboundaries in $W^{9}(SM)$.

\begin{proposition}\label{200}
Let $\mu \in spec(\Box)$, $T > 0$ and $k\in \mathbb{Z}$.  If $f\in W^{5}(\mathcal{H}_\mu)$ has a transfer function $u\in \mathcal{H}_\mu$, then $\hat{\delta}_{k/T}(f) = 0.$
\end{proposition}

The cases when $\mu \geq 1$ and $\mu \leq 0$ are similar, so we handle them together.

\begin{lemma}\label{990}
Let $\mu \leq 0$ or $\mu \geq 1$, and $T > 0$.  For all $k \in \Z$, if $f\in W^{5}(\mathcal{H}_\mu)$ has a transfer function $u\in \mathcal{H}_\mu$, then $\hat{\delta}_{k/T}(f) = 0.$
\end{lemma}
$\textit{Proof}:$
First suppose that $\mu \geq 1$, and let $u\in \mathcal{H}_\mu$ be such that $f = u\circ\phi_T - u.$  By extending the Fourier transform on $W^{1}(\mathcal{H}_\mu)$ as in definition $(\ref{224})$, we see that $\hat{f}$ is continuous.  Note that $\mathcal{H}_\mu$ takes the $L^2(\mathbb{R})$ norm, so $$\hat{f} = \mathcal{F}(u\circ \phi_T^U - u) = (e^{-2\pi i T \xi} - 1)\hat{u},$$ in $L^2(\mathbb{R})$.  Therefore, $\hat{u} = \frac{\hat{f}}{(e^{-2\pi i T \xi} - 1)}$ in $L^2(\mathbb{R})$.  Because $\hat{f}$ is continuous and $\hat{u}\in L^2(\mathbb{R})$, we conclude $\hat{f}(\frac{k}{T}) = 0,$ for all $k \in \Z$.  \\

When $\mu \leq 0$, again suppose $\hat{\delta}_{k/T}(f) \neq 0$, and recall the norm for the model $L^2(H, d\lambda_\nu)$ is $\| f \|_{L^2(H, d\lambda_\nu)} = \int_\mathbb{R}\int_\mathbb{R} |f(x + i y)|^2 y^{\nu - 1} dx dy.$  The same argument gives Lemma $\ref{990}$.  $\ \ \Box$\\

\begin{lemma}\label{1007}
Let $0 < \mu < 1$, $T > 0$ and $k\in \mathbb{Z}$.  If $f\in W^{5}(\mathcal{H}_\mu)$ has a transfer function $u\in \mathcal{H}_\mu$, then $\hat{\delta}_{k/T}(f) = 0.$
\end{lemma}
$\textit{Proof}:$
Suppose to the contrary that $\hat{\delta}_{k/T}(f) \neq 0$.  Let $K(x) = |x|^{-\nu}$ and $\psi \in \mathcal{S}(\mathbb{R}).$  Then $$\langle f, \psi\rangle_{\mathcal{H}_\mu} = \int_{\mathbb{R}^2} \frac{f(x) \overline{\psi(y)}}{|x - y|^{1 - \nu}} dx dy$$ \begin{equation}\label{1004} = \int_\mathbb{R} \left(\int_\mathbb{R} f(x) \frac{1}{|y - x|^{1 - \nu}} dx\right) \overline{\psi(y)} dy = \langle f\ast K, \psi\rangle_{L^2}.\end{equation}  Then define the linear functional $\ell_f$ by $\ell_f(\psi) := \langle f\ast K, \psi\rangle_{L^2(\R)},$ and similarly define $\ell_u$.  By Lemma $\ref{1095}$, \begin{equation}\label{3001}|\ell_f(\psi)| \leq \| f \|_{\mathcal{H}_\mu}\|\psi\|_{\mathcal{H}_\mu} \leq \| f \|_{\mathcal{H}_\mu}(\|\psi \|_{L^1(\mathbb{R})} + \| \psi \|_{L^\infty(\mathbb{R})}).\end{equation}  In particular, $\ell_f, \ell_u \in \mathcal{S}'(\mathbb{R}).$

By assumption $0 = \ell_{f - (u\circ\phi_T^U - u)},$ and by linearity it follows that $\ell_f = \ell_{u \circ\phi_T^U} - \ell_u.$  Observe that $\phi_T^U$ is unitary on $\mathcal{H}_\mu$, so $\ell_f = \phi_T^U \ell_u - \ell_u,$ and therefore \\ $\hat{\ell_f} = (e^{2\pi i \xi T} - 1) \hat{\ell_u},$ which means \begin{equation}\label{1005}\frac{\hat{\ell_f}}{(e^{2\pi i T \xi} - 1)} = \hat{\ell_u}\end{equation} in $\mathcal{S}'(\mathbb{R}).$

Now $\hat{\ell_f} = \hat{f}\hat{K}$ in $\mathcal{S}'(\mathbb{R})$, and we show that $\hat{\ell_f}$ is continuous away from 0.  Switching to circle coordinates, write $f\circ \tan(\theta) = \Phi(\theta) \cos^{1 + \nu}(\theta).$  Using Sobolev embedding followed by Lemma $\ref{554}$ gives $$\|\Phi\|_{C^0} \leq C_{\epsilon} \left(\sum_{k = -\infty}^{\infty} (1 + |k|)^{1 + \epsilon} |c_k|^2 \right)^{1/2}$$ $$\leq C_{\epsilon, SM}\| f \|_{W^{(1 + \Re\nu)/2 + \epsilon}(\mathcal{H}_\mu)} < \infty.$$ Then $$\| f \|_{L^1} \leq \|\Phi\|_{C^0} \int_{(-\pi/2, \pi/2)} \cos^{1 + \nu}(\theta) \frac{d\theta}{\cos^2(\theta)} < \infty,$$ so that $\hat{f}$ is continuous.  

Next, we show that for all $\xi \in \mathbb{R} - \{0\}$, $\hat{K}(\xi) \in \mathbb{C}$.  Notice $$|\hat{K}(\xi)| \leq \int_{|x| \leq 1} |x|^{-1 + \nu} dx + |\int_{|x| > 1} |x|^{-1 + \nu} e^{-2\pi i x \xi} d\xi|$$ 
Because $\nu \in (0, 1)$, it suffices to consider the integral with domain $\{|x| > 1\}$.  Using integration by parts, we have $$\int_{x > 1} x^{-1 + \nu} e^{-2\pi i x \xi} dx
 = \frac{-1 + \nu}{2\pi i \xi} \int_1^{\infty} x^{-2 + \nu} e^{-2\pi i x \xi} dx - \frac{1}{2\pi i \xi} \left(x^{-1 + \nu} e^{-2\pi i x \xi}\right)|_{1}^\infty \in \mathbb{C}.$$  
Moreover, if $\hat{K} \equiv 0$ on $\mathbb{R} - \{0\}$, then $K \equiv 0,$ which cannot be, so there exists some $\xi_0 \in \mathbb{R}$ such that $\hat{K}(\xi_0) \in \mathbb{C}^{\times}.$ 
 
 Let $\beta \in \mathbb{R} - \{0\}$, and notice that $K$ is homogeneous, so $$\hat{K}(\beta \xi_0) = \int_{-\infty}^{\infty} |x|^{-1 + \nu} e^{-2\pi i x \beta \xi_0} dx$$
  $$ = \text{sgn}(\beta) \int_{-\infty}^{\infty} |\frac{x}{\beta}|^{-1 + \nu} e^{- 2\pi i x \xi_0} \frac{dx}{\beta} = \text{sgn}(\beta) \frac{|\beta|^{1 - \nu}}{\beta} \hat{K}(\xi_0) = |\beta|^{-\nu} \hat{K}(\xi_0).$$  
 
By assumption $\hat{\delta}_{k/T}(f) \neq 0$, so if $k \neq 0$, then $\hat{f} \hat{K}$ is a continuous, nonzero function in a neighborhood $\mathcal{N}_k$ containing $k$.  Let $\psi \in C_c^{\infty}(\mathbb{R})$ be a nonzero bump function supported on $\mathcal{N}_k$ and  satisfying $\psi(k) \neq 0.$  Then $$\frac{\hat{\ell_f}}{(e^{2\pi i T \xi} - 1)}(\psi) = \hat{\ell_f}(\frac{\psi}{(e^{2\pi i T \xi} - 1)}) = \langle \hat{f} \hat{K}, \frac{\psi}{(e^{2\pi i T \xi} - 1)})\rangle_{L^2(\mathbb{R})}$$ \begin{equation}\label{1006} = \int_{\mathcal{N}_k} \frac{\hat{f}(\xi) \hat{K}(\xi)}{(e^{2\pi i T \xi} - 1)} \psi(\xi) d\xi = \infty.\end{equation}  If $k = 0$, we again conclude $\frac{\hat{\ell_f}}{(e^{2\pi i T \xi} - 1)}(\psi) = \infty.$
 
 On the other hand, $(\ref{3001})$ shows $$\hat{\ell_u}(\psi) \leq \| u \|_{\mathcal{H}_\mu} \left(\|\psi\|_{L^1(\mathbb{R})} + \| \psi\|_{L^{\infty}(\mathbb{R})}\right) < \infty.$$  But given $(\ref{1005})$ and $(\ref{1006})$, this is a contradiction.   $\ \ \Box$
 \\ \\
 $\textit{Proof of Proposition }\ref{200}:$  This follows immediately from Lemmas $\ref{990}$ and $\ref{1007}$.\\ 

Let $\mu > 0$ and $f\in W^{5}(\mathcal{H}_\mu)$.  The operator $A_T$ defined in $(\ref{320})$ represents as $$A_T f(x) = -\int_0^{-T} f(x + t) dt \text{ and } A_{T} f(z) = -\int_0^{-T} f(x + t + i y) dt$$ in the $\mathcal{H}_\mu$ model ($\mu > 0$) and the $L^2(H, d\lambda_\nu)$ model ($\mu \leq 0)$, respectively. 
We will use Proposition $\ref{333}$ to prove the distributions $\delta^{(0)}$ obstructs the existence of $\mathcal{H}_\mu$ solutions for $\mu \geq 0$, and we re-state it here for the convenience of the reader.   \\ \\
$\textit{Let }\triangle_M\textit{ have a spectral gap, and let } s > 1\text{ and }f\in W^{s}(SM)$.  $\textit{Then }$ $\textit{there exists }u\in L^2(SM)\textit{ such that }$ $$\mathcal{L}_U u = f \textit{ if and only if }u\circ \phi_T^U - u = A_T f.$$  

We defer its proof to Appendix A, Section $\ref{3022}$.  Proposition $\ref{333}$ implies $A_T$ maps coboundaries for the horocycle flow to coboundaries for the horocycle map.  Using Theorem $\ref{93}$, we find every smooth enough coboundary for the horocycle map arises this way.

\begin{lemma}\label{2070}
Let $T > 0$ and $f\in W^{9}(SM)$, and suppose there exists $u\in L^2(SM)$ such that \begin{equation}\label{202}f = u\circ\phi_T^U - u.\end{equation}  Then there exists $g\in W^{4/3}(SM)$ such that $A_T g = f.$
\end{lemma}
$\textit{Proof}:$
By $(\ref{304})$, write $f = \int_{\oplus \mu} f_\mu d\beta(\mu)$.  Now fix an irreducible component $\mathcal{H}_\mu$, and write $f = f_\mu$ for convenience.  By Lemma 6.3 of \cite{8} we have $Uf \in W^{8}(\mathcal{H}_\mu)$, and by flow invariance $U\delta^{(0)} = 0$.  Proposition $\ref{200}$ shows $f \in Ann(\{\hat{\delta}_k\}_{k = -\infty}^{\infty})$, so $$\hat{\delta}_{k/T}(Uf) = 2\pi i \frac{k}{T} \hat{\delta}_{k/T}(f) = 0.$$  Then Theorem $\ref{93}$ shows there exists $g\in W^{4/3}(\mathcal{H}_\mu)$ such that $Uf = g\circ \phi_T^U - g.$

If $\mu > 0$, then for every $M\in \mathbb{R}^-$, $$f(x) - f(M) = \int_{M}^x Uf(t) dt$$ $$ = \int_{M}^x [g(t - T) - g(t)] d = \int_{x}^{x - T} g(t) dt  - \int_{M}^{M - T} g(t) dt.$$  Write $g(\theta) = \Phi(\theta)\cos^{1 + \nu}(\theta)$.  Then by Sobolev embedding, $\|\Phi\|_{C^0(\mathbb{R})} 
\leq \| g \|_1$, so that for all $x\in \mathbb{R}$, $|g(x)| \leq \frac{\| g \|_1}{\sqrt{1 + x^2}}.$  Hence, $\lim_{M\rightarrow -\infty} \int_{M - T}^{M} g(t) dt = 0,$ and for the same reason, $\lim_{M\rightarrow-\infty} f(M) = 0,$ which means $$f(x) = \int_0^{-T} g(x + t) dt = A_T(-g)(x).$$   

If $\mu < 0$, we again see that for every $M \in \mathbb{R}^-$, $$f(x + i y) - f(M + i y) = \int_{x}^{x - T} g(t + i y) dt  - \int_{M}^{M - T} g(t + i y) dt.$$  Now Lemma $\ref{447}$ shows $|g(x + i y)| \leq \| g \|_1 \frac{1}{(1 + |x + i y|)^2},$ so $f = A_T(-g),$ as before.  $\ \ \Box$\\

\begin{proposition}\label{190} 
Let $\mu \geq 0$, $T > 0$, and $f\in W^{9}(\mathcal{H}_\mu)$, and suppose there exists $u\in \mathcal{H}_\mu$ such that $$f = u\circ\phi_T^U - u.$$ Then $\delta^{(0)}(f) = 0.$
\end{proposition}
$\textit{Proof}:$
By Lemma $\ref{2070}$, there exists $g\in W^{4/3}(\mathcal{H}_\mu)$ such that $f = A_T g$, and now Proposition $\ref{333}$ implies $\mathcal{L}_U u = g$.  With this, Lemmas 4.7, 4.8 and 4.9 of \cite{2}, show $\delta^{(0)}(g) = 0.$  Moreover, the flow invariance of $\delta^{(0)}$ implies that for all $t\in [0, 1]$, $\delta^{(0)}(g\circ \phi_t^U)$ $ = \delta^{(0)}(g) = 0,$ so $$0 = \int_0^{T}\delta^{(0}(-g \circ \phi_t^U) dt = \delta^{(0)}(-\int_0^{-T} g\circ \phi_t^U dt) = \delta^{(0)}(f). \ \ \Box$$\\

\begin{proposition}\label{617}
Let $\mu < 0$ and $r > 1$.  Then there exists $f\in W^r(H, d\lambda_\nu)$ such that $\delta^{(0)}(f) \neq 0$ and $f = u \circ\phi_T^U - u$, for some $u \in L^2(H, d\lambda_\nu)$.
\end{proposition}
$\textit{Proof}:$
We safely restrict ourselves to the holomorphic discrete series, so $\mu < 0$ implies $n = \frac{\nu + 1}{2} \geq 2$.  Let $f\in W^r(H, d\lambda_\nu)$ be such that $\delta^{(0)}(f) \neq 0.$  Then Lemma 4.5 of \cite{2} shows there is a solution $u\in L^2(H, d\lambda_\nu)$ such that $\mathcal{L}_U u = f$.  So Proposition $\ref{333}$ implies $$A_T f = u \circ\phi_T^U - u,$$ and notice that $$\|A_T f \|_r \leq \int_0^{T}\| f \circ\phi_t^U\|_{r} dt \leq C_{r, T} \| f \|_r,$$ by Minkowski's inequality and the commutation relations.  Finally, notice that $\delta^{(0)}$ is flow invariant, so $$\delta^{(0)}(A_T f) = -\delta^{(0)}\left(\int_0^{-T} f(\cdot + t) dt\right) = -\int_0^{-T} \delta^{(0)}(f(\cdot + t)) dt \neq 0. \ \ \Box \\$$

$\textit{Proof of Theorem } \ref{302}$:  Combining Propositions $\ref{200}$, $\ref{190}$ and $\ref{617}$, we conclude. $\ \ \Box$ \\

\section{Equidistribution of Horocycle Maps}

In this section we prove the equidistribution result listed as Theorem $\ref{152}$.  We focus on the time-1 horocycle map for simplicity, but the same argument works for the time-$T$ horocycle map.  Let $SM$ be compact, and note the Laplacian has only pure point spectrum, so we have the correspondence $$\mathcal{D}_n := (Q_\mu)^*\hat{\delta}_n, \ \  \mathcal{D}^0 := (Q_\mu)^* \delta^{(0)}$$ between the invariant distributions in $\langle\{\mathcal{D}_n\}_{n\in \mathbb{Z}}\cup \{\mathcal{D}^0\}\rangle = \mathcal{I}(\mathcal{K}_\mu)$ and those in $\langle\{\hat{\delta}_n\}_{n\in \mathbb{Z}}\cup \{\delta^{(0)}\}\rangle = \mathcal{I}(\mathcal{H}_\mu)$.

We assume $\int_{SM} f d vol = 0$ and prove that for all $(x_0, N)\in SM \times \mathbb{Z}^+$, \begin{equation}\label{1054}\frac{1}{N}\sum_{k = 0}^{N - 1} f(\phi_k^Ux_0) = \bigoplus_{\mu \in spec(\Box)} \sum_{\mathcal{D}\in \mathcal{I}^s(\mathcal{K}_\mu)} c_\mathcal{D}(x_0, N, s)\mathcal{D}(f) \oplus \mathcal{R}(x_0, N, s)(f),\end{equation} where $\mathcal{R}(x_0, N, s) \in (\mathcal{I}(SM))^\bot$ is the remainder distribution and satisfies $\| \mathcal{R}(x_0, N, s)\|_{-s} \leq \frac{C_s}{N}$, $c_\mathcal{D}(x_0, N, s)$ are the coefficients at the invariant distributions $\mathcal{D}$ and satisfy $|c_\mathcal{D}(x_0, N, s)| \leq C_s N^{-\mathcal{S}_\mathcal{D}}$ and $\mathcal{S}_\mathcal{D}$ is as in Theorem $\ref{152}$. \\
\subsection{Remainder distribution}

\begin{proposition}\label{156}
Let $s \geq 6$.  Then there exists a constant $C_s > 0$ such that for all $(x_0, N)\in SM\times \mathbb{Z}^+$, $$\|\mathcal{R}(x_0, N, s)\|_{W^{-s}(SM)}\leq \frac{C_s}{N}.$$
\end{proposition}
$\textit{Proof}:$
Let $f\in W^s(SM)$.  Because $\mathcal{R}(x_0, N, s)\in (\mathcal{I}^s(SM))^{\bot}$, we can write $f = f_\mathcal{I}\oplus f_\mathcal{C},$ where $f_\mathcal{I}\in (Ann(\mathcal{I}^s(SM)))^{\bot} \subset Ker(\mathcal{R}(x_0, N, s))$, and $f_\mathcal{C}\in Ann(\mathcal{I}^s(SM))$ is the coboundary component.  Then by the splitting in $(\ref{1054})$, \begin{equation}\label{154}|\mathcal{R}(x_0, N, s) f|  
 = |\mathcal{R}(x, N, s) f_\mathcal{C}| = |\frac{1}{N}\sum_{n = 0}^{N -1} \delta_{\phi_n^U(x_0)} f_{\mathcal{C}}|.\end{equation}  Theorem $\ref{93}$ and Sobolev embedding show there exists $\frac{1}{2} < r < s, C_{s} > 0$ and a (unique) $L^2(SM)$ solution $u\in W^{r}(SM)$ to the cohomological equation $u\circ \phi_T - u = f_\mathcal{C}$ satisfying $\|u\|_{C^0(SM)} \leq C_r \|u\|_r \leq C_{s} \| f_\mathcal{C} \|_s.$ Now $(\ref{154})$ becomes a telescoping sequence, and using this estimate, we conclude. $\ \ \Box$
\subsection{Invariant distributions}

For each $n \in \mathbb{Z}$, define $c_n := c_{\mathcal{D}_n}$, and write $d_0 := c_{\mathcal{D}^0},$ where $c_{\mathcal{D}_n}$ and $c_{\mathcal{D}^0}$ are the coefficients given in the asymptotic expansion ($\ref{1054}$).  The following lemma allows us to control the coefficients at the invariant distributions in terms of the horocycle flow and the "twisted" horocycle flow.  

\begin{lemma}\label{1080}
Let $\mu \in spec(\Box)$ and $s \geq 6$.  Then for all $\tau \in \mathbb{Z} - \{0\}$, $$c_\tau(x_0, N, s) \mathcal{D}_\tau = \frac{1}{N}\int_0^N e^{2\pi i \tau t} (\phi_t^U(x_0))^* dt - \int_0^1 e^{2\pi i \tau t} \phi_{-t}^U \mathcal{R} dt,$$ $$c_0(x_0, N, s) \mathcal{D}_0 = \frac{1}{N}\int_0^N (\phi_t^U(x_0))^* dt  - d_0(x_0, N, s) \mathcal{D}^0 - \int_0^1 \phi_{-t}^U \mathcal{R} dt,$$ and $$d_0(x_0, N, s) \mathcal{D}^0 = \frac{1}{N}\int_0^N (\phi_t^U(x_0))^* dt  - c_0(x_0, N, s) \mathcal{D}_0 - \int_0^1 \phi_{-t}^U \mathcal{R} dt$$ as distributions in $\mathcal{E}'(\mathcal{K}_\mu).$  
\end{lemma}
The proof of Lemma $\ref{1080}$ is contained in Appendix B.\\ \\

At this point, we use the (sharp) estimate of $\frac{1}{N}\int_0^N (\phi_t^U(x_0))^* dt$ given in Theorem 1.5 of \cite{2} to estimate the coefficients at the flow invariant distributions.
\begin{proposition}\label{950}
Let $\mu\in spec(\Box)$, $s \geq 6$ and the ergodic sum $(\ref{1054}) \in W^{-s}(SM)$.  Then there exists a constant $C_s > 0$ such that for all ($x_0, N)\in SM\times \mathbb{N}$, $$|c_0(x_0, N, s)| \leq C_s N^{-(1 - \Re\sqrt{1 - \mu})/2} \log^+(N),$$ and $$|d_0(x_0, N, s)| \leq C_s \bigg\{\begin{array}{ll} N^{-(1 + \Re\sqrt{1 - \mu})/2} \log^+(N) \text{ if } \mu > 0 \\ \frac{\log^+(N)}{N} \text{ if } \mu \leq 0.\end{array}$$  
\end{proposition}
$\textit{Proof}:$
We estimate $|c_0(x_0, N, s)|$.  Let $f_0 \in C^{\infty}(\mathcal{K}_\mu)$ be such that $\mathcal{D}_0(f_0) = 1$ and $\mathcal{D}^0(f_0) = 0.$  Let $\pi : W^1(\mathcal{K}_\mu) \rightarrow \left(Ann(\mathcal{I}(\mathcal{K}_\mu))\right)^{\bot}$ be orthogonal projection.  Let $g \in C^{\infty}(\mathcal{K}_\mu)\cap Ann(\mathcal{I}(\mathcal{K}_\mu))$, so there is a transfer function $u \in W^1(\mathcal{K}_\mu)$ corresponding to the coboundary $g$ such that for any $t\in [0, 1]$, $g\circ \phi_{-t}^U = u\circ \phi_{1 - t}^U + u\circ \phi_{-t}^U$.  Therefore, $g \circ\phi_{-t}^U \in Ann(\mathcal{I}(\mathcal{K}_\mu)).$ Hence, $$\langle \pi f_0\circ \phi_t^U, g\rangle_{\mathcal{K}_\mu} = \langle \pi f_0, g\circ \phi_{-t}^U\rangle_{\mathcal{K}_\mu} = 0.$$  Then for all $t\in [0, 1]$, $\pi f_0 \circ \phi_t^U \in \left(Ann(\mathcal{I}(\mathcal{K}_\mu))\right)^\bot$, so $\mathcal{R}(\pi f\circ\phi_t^U) = 0$. 
 
Then Lemma $\ref{1080}$ gives $$|c_0(x_0, N, s)| \leq |\frac{1}{N}\int_0^N \pi f_0(\phi_t^U(x_0)) dt| + |\int_0^1\phi_t^U(\mathcal{R})(x_0, N, s)(\pi f_0) dt| $$ \begin{equation} \label{1051}\leq |\frac{1}{N}\int_0^N \pi f_0(\phi_t^U(x_0)) dt|.\end{equation}  Finally, Theorem 1.5 of \cite{2} gives the estimate. 

The estimate for the coefficient $d_0(x_0, N, s)$ follows in the same manor. $ \ \ \Box$\\ \\

To estimate the coefficients at the other invariant distributions, recall that $\alpha(\mu_0) := \frac{(1 - \Re\sqrt{1 - \mu_0})^2}{8(3 - \Re\sqrt{1 - \mu_0})},$ where $\mu_0 > 0$ is again the spectral gap of $\triangle_{M}$.

\begin{proposition}\label{916}
Let $\mu \in spec(\Box)$, $\tau \in \mathbb{Z},$ $s \geq 6,$ $\epsilon > 0$ and the ergodic sum $(\ref{1054}) \in W^{-s}(SM).$  Then $$|c_\tau(x_0, N, s)| \leq C_\epsilon \tau^{5/2 + \epsilon} N^{-\alpha(\mu_0)}.$$
\end{proposition} 

 Define $$\gamma_{N, \tau} = \frac{1}{N}\int_0^N e^{2\pi i \tau t} (\phi_t^U(x_0))^* dt.$$  In light of Lemma $\ref{1080}$, Proposition $\ref{916}$ will follow once we estimate $\gamma_{N, \tau}$, which is given to us by a recent result of Venkatesh (Lemma 3.1 of \cite{23}).  For any natural number $k \geq 0$, let $W^{k, \infty}(SM)$ be the set of $L^2(SM)$ functions satisfying $$\| f \|_{W^{k, \infty}(SM)} = \| (1 + \triangle)^{k/2} f\|_{L^\infty(SM)} < \infty.$$

\begin{lemma}\label{925}[Venkatesh] 
Let $\int_{SM} f dvol = 0$, $\epsilon > 0$ and $\mu_0 > 0.$  Then there exists a constant $C > 0$ such that for all $(x_0, N) \in SM \times \mathbb{N}$ and $\tau \in \mathbb{R},$ we have $$|\gamma_{N, \tau}(f)| \leq C \| f \|_{W^{1, \infty}(SM)} N^{-\alpha(\mu_0)}.\\$$
\end{lemma}
$\textit{Proof of Proposition }\ref{916}:$
First let $\mu \in spec(\Box)$.  Lemma $\ref{715}$ gives a function $f_\tau \in C^\infty(\mathcal{K}_\mu)$ such that for all $\mathcal{D} \in \mathcal{I}(\mathcal{K}_\mu)$, $$\mathcal{D}(f_\tau) = \Big\{\begin{array}{rr}1 \text{ if }\mathcal{D} = \mathcal{D}_\tau \\0 \text{ otherwise}.\end{array}.$$  Let $\pi : W^1(\mathcal{K}_\mu) \rightarrow \left(Ann(\mathcal{I}(\mathcal{K}_\mu))\right)^{\bot}$ be orthogonal projection as in Proposition $\ref{950}$.

Then Lemma $\ref{1080}$ together with Lemma $\ref{925}$ gives $$|c_\tau(x, N, s)| \leq C_\epsilon \|f_\tau\|_{W^{1, \infty}(\mathcal{K}_\mu)} N^{-\alpha(\mu_0)} + |\int_0^1 e^{2\pi i \tau t} \phi_{-t}^U\mathcal{R}( \pi f_\tau) dt|$$ $$ \leq C \| f_\tau \|_{W^{1, \infty}(\mathcal{K}_\mu)} N^{-\alpha(\mu_0)}.$$  Finally, Sobolev embedding and the estimates $(\ref{2101})$ and $(\ref{2100})$ in the proof of Lemma $\ref{715}$ prove that for any $\epsilon > 0$, $$\| f_\tau \|_{W^{1, \infty}(\mathcal{K}_\mu)} \leq C_\epsilon \| f_\tau \|_{W^{5/2 + \epsilon}(\mathcal{K}_\mu)} \leq C_\epsilon \tau^{5/2 + \epsilon}. \ \ \Box$$

\subsection{Proof of Theorem $\ref{152}$}

Propositions $\ref{156}$, $\ref{950}$, and $\ref{916}$ give an upper bound for the rate of decay of the remainder distribution and all invariant distributions.  Now we need conditions showing when the series in $(\ref{1054})$ converges.  We begin with a lemma, whose proof is deferred to Appendix B.

\begin{lemma}\label{787}
Let $\mu\in spec(\Box), s \geq 2$ and $f\in W^s(\mathcal{H}_\mu).$  Then for all $\xi\in \mathbb{R}$, $$|\hat{f}(\xi)| \leq C_s \|f\|_{W^{3s + 2}(\mathcal{H}_\mu)} (1 + |\xi|)^{-s}.$$  
\end{lemma}

$\textit{Proof of Theorem }\ref{152}:$  Recall that $\alpha(\mu_0) := \frac{(1 - \sqrt{1 - \mu_0})^2}{8(3 - \sqrt{1 - \mu_0})}.$  We estimate one irreducible component at a time, so let $\mu  > 0$ and $f \in W^{s}(\mathcal{K}_\mu)$, where $s \geq 14$.   Then by Lemma $\ref{787}$ and Propositions $\ref{156}$, $\ref{950}$, and $\ref{916}$, we have  
$$|\frac{1}{N}\sum_{n = 0}^{N - 1} f(\phi_n^Ux_0)| = |\left(\sum_{k \in \mathbb{Z}}c_k(x_0, N, s)\mathcal{D}_k(f)  + d_0(x, N, s) \mathcal{D}^0(f) \right) \oplus \frac{\mathcal{R}(x, N, s)(f)}{N}|$$ $$\leq C_{s} \| f \|_{s} \sum_{k \in \mathbb{Z} - \{0\}} |k|^{5/2 + 1/4 } |k|^{-4} N^{-\alpha(\mu_0)} \log^+(N)$$ $$+ C_s \| f \|_s \left(N^{(1 - \nu)/2}\log^+(N) + N^{-(1 + \nu)/2} \log^+(N) + \frac{1}{N}\right).$$  One estimates in the same way when $\mu \leq 0$.  The series converges absolutely with constant $C_s$ independent of $\mu$.  Then Theorem $\ref{152}$ follows by gluing the series together from each irreducible component. $ \ \ \Box$

\appendix
\section{}

\subsection{Claim $\ref{133}$ and proof of Lemma $\ref{554}$}

The models for the principal, complementary and discrete series are discussed in Section 2.  For $\mu > 0$, a calculation in the circle model shows $\cos^{1 + \nu}(\theta) \in Ker (\Theta)$.  For all $k \in \Z$, define $$u_k := e^{2 i k \theta} \cos^{1 + \nu}(\theta).$$  For $\mu \leq 0$, we do calculations in the upper half-plane model $L^2(H, d\lambda_\nu)$.  A formal calculation gives $u_0 = \left(\frac{z - i}{z + i}\right)^{-n} \left(\frac{1}{z + i}\right)^{\nu + 1} \in Ker \Theta,$ where $n = \frac{\nu + 1}{2}.$  In this case, for all integers $k \geq n$, define $$u_k := \left(\frac{z - i}{z + i}\right)^{k - n} \left(\frac{1}{z + i}\right)^{\nu + 1}.$$

Let $j_\nu = \Big\{\begin{array}{cc} - \infty &\text{ if } \mu > 0\\  n &\text{ if } \mu \leq 0. \end{array}$  
\begin{claim}\label{133}
For any $\mu \in spec(\Box)$ and for all $k \geq j_\nu$, we have
$$\begin{array}{lll} i) \ (X + i Y) u_k = -(1 + \nu + 2k) u_{k + 1} \text{ and } (X - i Y) u_k = [-(1 + \nu) + 2k] u_{k - 1}.\\
ii) \ \ - i \Theta(u_k) = 2k u_k.\\
iii) \ \ \ \Box u_k = (1 - \nu^2) u_k \text{ and } \triangle u_k = (1 - \nu^2 + 8 k^2) u_k.\end{array}$$ 
Additionally, $\{u_k\}_{k \geq j_\nu}\subset C^{\infty}(\mathcal{H}_\mu)$ (resp. $C^\infty(H, d\lambda_\nu)$) is an orthogonal basis for $\mathcal{H}_\mu$ (resp. $L^2(H, d\lambda_\nu)$). $ \ \ \Box$
\end{claim}




$\textit{Proof of Lemma }\ref{554}:$

Part $i):$  
By Claim $\ref{133}$, $\{u_k = e^{2i k \theta} \cos^{1 + \nu}(\theta)\}_{k \in \mathbb{Z}}$ is an orthonormal basis for $\mathcal{H}_\mu$.  If $\mu \geq 1$, we may write $\nu = i s$ for $s \in \mathbb{R}$ and get $$\|u_k\|^2 = \int_{-\pi/2}^{\pi/2} |e^{2 i k \theta} \cos^{1 + is}(\theta)|^2 \frac{d\theta}{\cos^{2}(\theta)} = \pi.$$  

When $0 < \mu < 1$, we see $$\|u_0\|_{\mathcal{H}_\mu}^2 = \int_{\mathbb{R}^2} \frac{(1 + x^2)^{-(1 + \nu)/2} (1 + y^2)^{-(1 + \nu)/2}}{|x - y|^{1 - \nu}} dx dy < \infty,$$ by splitting the integral into parts where $|x - y| \leq 1$ and $|x - y| \geq 1$.  Because there are only finitely many eigenvalues of $\Box$ in $(0, 1)$, we conclude there is a constant $C_{SM} > 0$ such that $C_{SM}^{-1} \leq \|u_0\| \leq C_{SM}.$

The basis $\{\tilde{u}_k\}$ in Flaminio-Forni \cite{2} is constructed using a vector $\tilde{u}_0 \in Ker(\Theta)$ normalized so that $\|\tilde{u}_0\|^2 = 1$, and then generating the rest of the elements from $\tilde{u}_0$ by the creation and annihilation operators.  Analogously to Claim $\ref{133}$, Formula (24) of \cite{2} gives $\|\eta_{\pm} \tilde{u}_k\| = \|(1 + \nu \pm 2k) \tilde{u}_{k \pm 1}\|.$  Hence, for all $k \geq 0$, $$\|u_{k + 1}\| = \| \frac{\eta_+ u_k}{1 + \nu + 2k}\| = \|\Pi_{j = 0}^k \frac{1}{1 + \nu + 2j} (\eta_+)^{k + 1} u_0\|$$ $$ = \|\tilde{C}_\nu \Pi_{j = 0}^k \frac{1}{1 + \nu + 2j} (\eta_+)^{k + 1} \tilde{u}_0\| = \tilde{C}_\nu \|\tilde{u}_{k + 1}\|,$$ where $C_{SM}^{-1} \leq \tilde{C}_\nu \leq C_{SM}$.  Lemma 2.1 of \cite{2} gives that whenever there is a spectral gap $\mu_ 0 > 0$, there is a constant $C_{SM} > 0$ such that $$C_{SM}^{-1} (1 + |k|)^{-\nu} \leq \|\tilde{u}_k\|^2 \leq C_{SM} (1 + |k|)^{-\nu}.$$ 

Part $ii):$
The unit disc model $L^2(D, d\sigma_\nu)$ has the measure $d\sigma_\nu := 4\frac{(1 - |\xi|^2)^{\nu - 1}}{|\xi - 1|^{2(\nu + 1)}} du dv$, and we use the conformal map $\alpha(\xi) = -i \left(\frac{\xi + 1}{\xi - 1}\right).$  One checks that $u_n\circ \alpha(\xi) = \left(\frac{\xi - 1}{-2i}\right)^{\nu + 1}.$  Using polar coordinates, $$\| u_n\|_{L^2(H, d\lambda_\nu)}^2 = \| u_n\|_{L^2(D, d\sigma_\nu)}^2 = 4 \int_D \left(\frac{|\xi - 1|}{2}\right)^{2(\nu + 1)}\frac{(1 - |\xi|^2)^{\nu - 1}}{|\xi - 1|^{2(\nu + 1)}} du dv$$ $$ = \pi 4^{-\nu} \int_0^1 t^{\nu - 1} dt = \frac{\pi}{\nu 4^{\nu}}.$$

Let $\{\tilde{u}_k\}_{k = n}^\infty$ be the basis given in \cite{2}.  Lemma 2.1 of \cite{2} gives that for all $k \geq n$, $$\|\tilde{u}_k\|^2 = \tilde{\Pi}_{\nu, k} = \frac{(k - n)! \nu!}{(k + n - 1)!}.$$  In particular, $\|\tilde{u}_n\|^2 = \tilde{\Pi}_{\nu, n} = 1,$ so $\| u_n \| = \frac{\sqrt{\pi}}{\sqrt{\nu} 2^{\nu}}\| \tilde{u}_n\|.$  

We generate the other basis vectors from our creation operator.  For all $k > n$, $$u_k = \frac{1}{\nu  + 2k - 1} \eta_+ u_{k - 1}\text{ and }\tilde{u}_k = \frac{1}{\nu + 2k - 1} \eta_+ \tilde{u}_{k - 1}.$$  By iterating we conclude $$\| u_k \| = \sqrt{\frac{\pi}{\nu}} 2^{-\nu} \| \tilde{u}_k \| = \sqrt{\frac{\pi}{\nu}} 2^{-\nu} \left(\frac{(k - n)! \nu!}{(k + n - 1)!}\right)^{1/2}.$$  With this, Lemma $\ref{554}\ ii)$ follows. $\ \ \Box$ \\   

\subsection{Proof of Theorem $\ref{322}$}
 
 The first step is to prove
  \begin{lemma}\label{733}
Let $\mu \leq 0$ and $r, s$ be integers such that $0 \leq r \leq s \leq \nu$.   If $f\in W^s(H, d\lambda_\nu),$ then for all $\xi\in D$, \begin{equation}\label{346}|(f \circ\alpha)^{(r)}(\xi)| \leq C_{r}\nu^{r + 1/2} \|f\|_s \sum_{j = 0}^r \left(\sum_{k = r - j}^\infty\frac{(k + \nu)!}{k! \nu!}(k + \nu)^{-2s + 2r} |\xi|^{2k}\right)^{1/2} |\xi|^{-r} |1 - \xi|^{(\nu + 1) - r}\end{equation}
\end{lemma}
$\textit{Proof}:$
We have $$f\circ\alpha(\xi) = (-2i)^{-(\nu + 1)} \sum_{k = 0}^{\infty} c_{k + n} \xi^k (1 - \xi)^{\nu + 1}.$$  Differentiating,
\begin{equation}\label{460} |f\circ\alpha^{(r)}(\xi)| \leq \frac{C_r}{2^{\nu}} \nu^r \sum_{j = 0}^r\left(\sum_{k = r - j}^\infty |c_{k + n}| \frac{k!}{(k - r + j)!}|\xi|^{k} \right) |\xi|^{-r + j} |1 - \xi|^{\nu + 1 - j}.\end{equation} 
Formula $(\ref{4015})$ gives the value of $\|u_{k + n}\|_s$.  Now multiplying and dividing by $\|u_{k + n}\|_s$ inside the sum and applying Cauchy-Schwarz, we conclude. $ \ \ \Box$ 

Now let $$B_T := \{z\in H : |z - i| < T/3\} \text{ and } B_T^{c, 0} := int(H - B_T).$$

 \begin{lemma}\label{447}
Let $\mu \leq 0$, and let $r, s$ be integers such that $0 \leq r \leq s \leq \nu$.  Also let $z\in B_T^{c, 0}$ and $f\in W^s(H, d\lambda_\nu)$.  If $\nu/2 + r < s$, then $$|f^{(r)}(z)| \leq C_{r}\left(\frac{1 + T}{T}\right)^r \frac{\nu^{r + 1/2}}{\sqrt{\nu!}} \| f \|_{s} \left( \frac{1}{\Im z}\right)^{1/2} |\frac{1}{1 + |z|^2 + 2\Im(z)}|^{(\nu + r)/2},$$ and if $\nu/2 + r \geq s$, then $$|f^{(r)}(z)| \leq C_{r} \nu^{r + 1/2} \left(\frac{1 + T}{T}\right)^r \| f \|_{s} \left( \frac{1}{\Im z} \right)^{\nu/2 - s + r + 1/2} \left( \frac{1}{1 + |z|^2 + 2\Im(z)}\right)^{ s - r/2}.$$  

Moreover, if $\mu = 0, \epsilon > 0$ and $f \in W^{1 + \epsilon}(\mathcal{H}_\mu)$, then $$|f(z)| \leq C_{\epsilon} \| f \|_{1 + \epsilon} \frac{1}{(1 + |z|)^2}.$$
\end{lemma}
$\textit{Proof}:$
By $(\ref{350})$ and Lemma $\ref{733}$ we have $$|f^{(r)}(z)| = U^r (f\circ\alpha)(\xi) \leq C_r \sum_{j = 1}^{r}|\xi - 1|^{r + j} |f^{(j)}(\xi)|$$ \begin{equation}\label{449} \leq C_{r}\frac{\nu^{r + 1/2}}{\sqrt{\nu!}} \|f\|_s \left(\sum_{k = r - j}^\infty \frac{(k + \nu)!}{k!} (k + \nu)^{-2s + 2r} |\xi|^{2k}\right)^{1/2} |\xi|^{-r} |1 - \xi|^{\nu + 1 + r}\end{equation} 

For the following Case 1 and Case 2, let $q = |\xi|^2$.  
\\ $\textit{Case 1}: \nu/2 + r < s$.  Then $\nu - 2s + 2r < 0,$ which means $\frac{(k + \nu)!}{k!} (k + \nu)^{-2s + 2r} \leq 1,$ and therefore $$\left(\sum_{k = 0}^{\infty} \frac{(k + \nu)!}{k!} (k + \nu)^{-2s + 2r} |\xi|^{2k} \right)^{1/2} \leq \left(\sum_{k = 0}^{\infty} q^{k} \right)^{1/2} = \left(\frac{1}{1 - q}\right)^{1/2}.$$ by the change of variable $\alpha:\xi\rightarrow -i \left(\frac{\xi + 1}{\xi - 1}\right),$ we have that \begin{equation}\label{614}q =  |\frac{x + i(y - 1)}{x + i(y + 1)}|^2 = 1 - \frac{4\Im z}{1 + |z|^2 + 2\Im z}\text{ and }|1 - \xi| = \left(\frac{4}{1 + |z|^2 + 2 \Im(z)}\right)^{1/2}.\end{equation}  Then noting that $z \in B_T^{c, 0}$ an elementary calculation concludes Case 1. \\
$\textit{Case 2}: \nu/2 + r \geq s.$  Then $\nu - 2s + 2r \geq 0$, so that $$\sum_{k = r - j}^\infty \frac{(k + \nu)!}{k!} (k + \nu)^{-2s + 2r} |\xi|^{2k} \leq \sum_{k = 0}^{\infty} \frac{(k + \nu - 2s + 2r)!}{k!} q^{k}$$ $$ = \frac{d^{\nu - 2s + r}}{dq^{\nu - 2s + 2r}}\left(\sum_{k = 0}^{\infty} q^k\right) = \frac{d^{\nu - 2s + 2r}}{dq^{\nu - 2s + 2 r}} \left(\frac{1}{1 - q}\right) = (\nu - 2s + 2 r)! \left(\frac{1}{1 - q}\right)^{\nu - 2s + 2 r + 1}.$$  Then combining this with $(\ref{449})$ and $(\ref{614})$ gives Case 2. 
 \\
 $\textit{Case 3}:$ If $\mu = 0$ and $f\in W^{1 + \epsilon}(\mathcal{H}_\mu)$, write $f \circ \alpha(\xi) = \Phi(\xi) (\xi - 1)^{\nu + 1}.$  Then Sobolev embedding gives $$|f(z)| = |f\circ\alpha(\xi)| = |\Phi(\xi)| |\xi - 1|^{\nu + 1} \leq \|\Phi\|_{C^0(D)} |\xi - 1|^{\nu + 1} $$ \begin{equation}\label{3030}\leq C_\epsilon \| f \|_{(1 + \nu)/2 + \epsilon} |\xi - 1|^{-(\nu + 1)} = C_\epsilon \| f \|_{1 + \epsilon} |\xi - 1|^{-2}. \ \ \Box \end{equation}  
 
  \begin{proposition}\label{748}
Let $y \in \mathbb{R}^+$, $s \geq 1$ and $f\in W^s(\mathcal{H}_\mu)$.  

If $\mu = 0$ and $s > 1$, then $f \in L^1(\cdot + i y)$.   If $\mu < 0$, and $r, s$ are integers such that $0 \leq r < s \leq \nu$, then $f^{(r)}(\cdot + i y) \in L^1(\mathbb{R}).$
 \end{proposition}
$\textit{Proof}:$  Because $f$ is holomorphic, it is bounded on compact sets, so it is bounded on $B_T\cap (-\infty, \infty) \times\{y\}.$  Then the proposition follows from Lemma $\ref{447}$. $\ \ \Box$ \\

$\textit{Proof of Lemma }\ref{627}:$
The statement that $\delta_{k/T, y} \in W^{-(1 + \epsilon)}(H, d\lambda_\nu)$ follows from Lemma $\ref{447}$ and $\phi_T^U$-invariance follows as in Lemma $\ref{901} \ \ \Box$.\\

$\textit{Proof of Lemma }\ref{767}:$
Say $y_1 > y_2$, and let $s > 1/2$.  Additionally, for all $n \in \mathbb{N}$, let $\Gamma_n$ be an oriented closed curve with sides $$\Gamma_n := \{[-n + i y_1, n + i y_1]\} \cup \{[n + i y_1, n + i y_2]\} \cup \{[n + i y_2, -n + i y_2]\} \cup \{[-n + i y_2, -n + i y_1]\}.$$  Let $f\in W^s(\mathcal{H}_\mu)$, and note that $f$ is holomorphic.  Then by Cauchy's theorem, $0 = \int_{\Gamma_n} f(z) e^{-2\pi i k/T z} dz.$  By Lemma $\ref{447}$, there is a constant $C_{\nu, y_1, y_2} > 0$ such that $$|\int_{[-n + i y_2, -n + i y_1]} f(z) e^{-2\pi i k/T z} dz| + |\int_{[n + i y_2, n + i y_1]} f(z) e^{-2\pi i k/T z} dz|$$ $$  \leq C_{\nu, y_1, y_2} \| f \|_s (1 + |n|)^{-2s}.$$ Letting $n \rightarrow \infty$, we conclude $\hat{\delta}_{k/T, y_1} = \hat{\delta}_{k/T, y_2}$.  

The second statement is proved the same way.  $\ \ \Box$\\

$\textit{Proof of Theorem} \ \ref{322}:$
For $\mu > 0$, formula (40) of \cite{2} and Lemma $\ref{901}$ and show that for all $\epsilon > 0$, we have $\delta^{(0)}, \hat{\delta}_{k/T} \in W^{-((1 + \nu)/2 + \epsilon)}(\mathcal{H}_\mu)$ are $\phi_T^U$-invariant distributions.  Moreover, the proof of Theorem $\ref{165}$ shows $\mathcal{I}(\mathcal{H}_\mu) = \langle \{\hat{\delta}_{k/T}\}_{k \in \Z} \cup \{\delta^{(0)}\}\rangle.$  Then unitary equivalence gives Theorem $\ref{322}$ for $\mu > 0$.  When $\mu \leq 0$, we use formula (40) of \cite{2} and Lemmas $\ref{767}$ and $\ref{627}$ for the regularity of $\delta^{(0)}$ and $\hat{\delta}_k$ and for the definition of $\hat{\delta}_k$, and we use Theorem $\ref{357}$ in place of Theorem $\ref{165}$. $ \ \ \Box$

\subsection{Proof of Proposition $\ref{929}$}\label{3020}

$\textit{Throughout this subsection},$ $\mu \leq 0$, $r, s\in \mathbb{N}_0, s \geq 4$ and $0 \leq r \leq \tilde{s}$, where $\tilde{s} := \lfloor \frac{s - 1}{2}\rfloor.$  Moreover, we use the conformal map $\alpha : D\rightarrow H : \xi \rightarrow -i\left(\frac{\xi + 1}{\xi - 1}\right) := z$ between $D$ and $H$, and we let $\{u_k\}_{k \geq n} \subset L^2(H, d\lambda_\nu)$ be the basis given in Lemma $\ref{554}$.

Note that if $0 \leq r \leq \tilde{s}$, then $\delta^{(r)} \in W^{-(r + (1 + \nu)/2 + \epsilon)}(H, d \lambda_\nu) \subset W^{-s}(H, d\lambda_\nu),$ because $r + \frac{1 + \nu}{2} + \epsilon \leq s$ when $0 < \epsilon \leq \frac{s - \nu}{2}.$ \\
\\
To begin the proof of Proposition $\ref{929}$, define the space $$\mathcal{P}_\nu(D) = \{\sum_{k = n}^M c_k u_k\circ\alpha  | M \geq n\text{ is an integer and }\{c_k\}_{k = n}^M \subset \mathbb{C}\},$$ and note $u_k\circ \alpha = \xi^{k - n} u_n(\xi)$.  For $f\in \mathcal{P}_\nu(D)$ write $f = \Phi \cdot u_n,$ where $\Phi(\xi) := \sum_{k = n}^M c_k \xi^{k - n}.$

\begin{lemma}\label{611}
If $f\in \mathcal{P}_\nu(D) \cap Ann(\{\delta^{(j)}\}_{j = 0}^{\tilde{s} - 1})$, then $\Phi^{(r)}(1) = 0.$
\end{lemma}
$\textit{Proof}:$
Because $f\in \mathcal{P}_\nu(D)$, we may differentiate $\Phi$ term by term.
Then \begin{equation}\label{29}\frac{d^r}{d\theta^r} \Phi(e^{2\pi i \theta})|_{\theta = 0} = \sum_{k \geq n} c_k (2\pi i k)^r = (-\pi)^r\delta^{(r)}(f) = 0.\end{equation} 

Now Taylor expand $\Phi$ about $\xi = 1$ and get $\Phi(\xi) = \sum_{k = 0}^{\infty} \beta_k (\xi - 1)^k,$ where $\{\beta_k\}_{k = 0}^{\infty}\subset\mathbb{C}$.  Then conclude $\beta_k = 0$ for all $0 \leq k \leq \tilde{s} - 1$ by induction.  $ \ \ \Box$

Let $j\in \mathbb{N}$ and given $t_{1}, \ldots, t_j\in \mathbb{R}$, let $\vec{t}_{j} := t_1 \cdots t_j \text{ and } \xi_{\vec{t}_j} := \vec{t}_j(\xi - 1) + 1.$

\begin{lemma}\label{35}
Let $f\in \mathcal{P}_\nu(D) \cap Ann(\{\delta^{(j)}\}_{j = 0}^{\tilde{s} - 1})$.  Then for all $\xi\in \mathcal{D}$, we have $$|\Phi^{(r)}(\xi)| \leq |\xi - 1|^{\tilde{s} - r} \int_{[0, 1]^{s - r}} |\Phi^{(\tilde{s})}(\vec{t}_{\tilde{s} - r} (\xi - 1) + 1)| dt_{\tilde{s} - r} \cdots dt_1.$$
\end{lemma}
$\textit{Proof}:$
Let $$g(t) = \Re\Phi^{(r)}(t(\xi - 1) + 1) + i \Im\Phi^{(r)}(t(\xi - 1) + 1).$$  Then the Fundamental Theorem of Calculus gives $$\Phi^{(r)}(\xi) = \int_0^1 g'(t_1) dt_1 = (\xi - 1) \int_0^1\Phi^{(r + 1)}(t_1(\xi - 1) + 1) dt_1,$$ and iterating gives Lemma $\ref{35}$.  $\ \ \Box$

Define $$\mathcal{P}_\nu(H) : = \{ f\in L^2(H, d\lambda_\nu) | f \circ \alpha \in \mathcal{P}_\nu(D)\}.$$

\begin{lemma}\label{472}
Let $f\in \mathcal{P}_\nu(H)\cap Ann(\{\delta^{(r)}\}_{r = 0}^{\tilde{s} - 1})$.  Then there is a constant $C_{r, s} > 0$ such that for all $z \in H$, $$|f^{(r)}(z)| \leq C_{r, s} \| f \|_s (1 + |z|)^{-(s/2 + \nu + r)}.$$ 
\end{lemma}
$\textit{Proof}:$  
Notice $f\circ \alpha \in \mathcal{P}_\nu(D)$, so $f\circ\alpha(\xi) = \Phi(\xi) (\xi - 1)^{\nu + 1}.$  Then using Lemma $\ref{35}$, we get $$(f\circ\alpha)^{(r)}(\xi) = \sum_{j = 0}^r (_j^r) \Phi^{(r - j)}(\xi) \frac{d^j}{d\xi^j} (\xi - 1)^{\nu + 1}$$
$$\leq C_s \nu^r \sum_{j = 0}^r  \int_{[0, 1]^{\tilde{s} - r + j}} |\Phi^{(\tilde{s})}(\xi_{\vec{t}_{\tilde{s} - r + j}})| d t_{\tilde{s} - r + j}\cdots d_{t_1} |\xi - 1|^{\tilde{s} + \nu + 1 - r}.$$ 

Because $D$ is convex, we know that $\xi_{\vec{t}_{\tilde{s} - j}} \in D$ for all $\vec{t}_{\tilde{s} - j}.$  Recall that $\nu < s$ and let $0 < \epsilon = \frac{s - \nu}{2}$.  For all $\xi_{\vec{t}_{\tilde{s} - j}} \in D$, $$|\Phi^{(\tilde{s})} (\xi_{\vec{t}_{\tilde{s} - j}})| \leq C_{\epsilon} \left(\sum_{k = -\infty}^{\infty} |c_k|^2 k^{2 \tilde{s} + 1 + \epsilon} \| u_k \|^2 \| u_k\|^{-2}\right)^{1/2}$$ $$\leq C_{\epsilon} \| f \|_{\tilde{s} + \frac{\nu + 1}{2} + \epsilon} \leq C_{s} \| f \|_{s}$$ 

One shows by induction that for $r \geq 1$, there are constants $\{c_j\}_{j = 1}^r\subset\mathbb{C}$ such that \begin{equation}\label{350}U^r (f\circ\alpha)(\xi) = \sum_{j = 1}^{r}c_j(\xi - 1)^{r + j} (f\circ\alpha)^{(j)}(\xi).\end{equation}  Then there exists $C_r > 0$ such that $$|U^r (f\circ \alpha)(\xi)| \leq C_{r} \sum_{j = 0}^r |\xi - 1|^{r + j} |(f\circ\alpha)^{(j)}(\xi)|$$ $$ \leq C_{r, s} \nu^r \| f \|_{s}\sum_{j = 0}^r |\xi - 1|^{\tilde{s} + \nu + 1 - j + j + r} \leq C_{r, s} \nu^r \| f \|_{s}\sum_{j = 0}^r |\xi - 1|^{(s/2 + \nu + r)}.$$ Then the Lemma follows from the change of variable given by the Mobius transformation $\alpha : \xi \rightarrow -i\frac{\xi + 1}{\xi - 1}$ (see $(\ref{614})$). $\ \ \Box$ \\

$\textit{Proof of Proposition }\ref{929}:$
Clearly, $\mathcal{P}_\nu(H)$ is dense in $W^s(H, d\lambda_\nu)$.  Then let $\eta > 0$ and $f_ \eta \in \mathcal{P}(H, d\lambda_\nu)\cap Ann(\{\delta^{(r)}\}_{r = 0}^{\tilde{s} - 1})$ satisfy $$\|f - f_ \eta\|_{W^{s}(H, d\lambda_\nu)} < \eta.$$  As $r < \frac{s - 1}{2}$ and $\nu < s$, take $0 < \epsilon = \frac{s - \nu}{2}$ and conclude that for all $z\in H$,
$$|\frac{d^r}{dz^r}(f - f_{\eta}) (z) | 
\leq \| \frac{d^r}{dz^r}(f - f_{\eta}) (z)\|_{(1 + \nu)/2 + \epsilon} \leq \|f - f_{\eta}\|_{s} < \eta,$$ where we use Nelson \cite{11} in the second inequality. 

Hence, $$|f^{(r)}(z)| \leq |\frac{d^r}{dz^r} (f - f_{\eta})(z)| + |f_{\eta}^{(r)}(z)|$$ $$ \leq \eta + C_{r, s} \cdot \nu^{r} (|z| + 1)^{-(s/2 + \nu + r)} \|f_{\eta}\|_{W^{s}(\mathcal{H}_\mu)}.\ \ \ \Box$$
 
 
\subsection{Proof of Proposition $\ref{333}$}\label{3022}

The left equality implies the right by the fundamental theorem of calculus, so it is the other direction that is of interest.  We consider each irreducible component individually, and the first step is the following lemma. 

\begin{lemma}\label{367}
Let $\mu \in spec(\Box)$.  If $\mu > 0$, then (distributional) $Ker(A_T) = \langle\{ \hat{\delta}_{n/T}\}_{n \in \mathbb{Z} - \{0\}}\rangle$ in the $\mathcal{H}_\mu$ model.  If $\mu \leq 0$, then (distributional) $Ker(A_T) = \langle\{ \hat{\delta}_{k/T}\}_{k \in \mathbb{Z}^+}\rangle$ in $L^2(H, d\lambda_\nu)$.  
\end{lemma}
 $\textit{Proof}:$
First let $\mu > 0$.  Let $\mathcal{D}\in Ker(A_T) \subset \mathcal{E}'(\mathcal{H}_\mu)$ and $h\in \mathcal{S}(\mathbb{R}) \subset C^{\infty}(\mathcal{H}_\mu).$  Recall that the Fourier transform is defined on $C^{\infty}(\mathcal{H}_\mu)$ by $(\ref{224})$, so $$0 = (A_T \mathcal{D}) h = -\mathcal{D}(A_{-T} h)$$ \begin{equation}\label{3007} = \hat{\mathcal{D}}(\int_0^{T} e^{2\pi i t \xi} dt \hat{h}(\xi)) = \hat{\mathcal{D}}(\frac{e^{2\pi i T \xi} - 1}{2\pi i \xi} \hat{h}(\xi)).\end{equation}  Hence, $supp(\hat{\mathcal{D}}) \subset \mathbb{Z} - \{0\}$, so $\langle \{\hat{\delta}_{k/T}\}_{k\in \mathbb{Z} - \{0\}}\rangle \subset Ker(A_T).$

For the other inclusion, $(\ref{3007})$ gives $Ker(A_T) \subset \langle \{\hat{\delta}_{k/T}^{(j)}\}_{k \in \mathbb{Z} - \{0\}, j \in \mathbb{N}_0}\rangle.$  Because distributions in $Ker(A_T)$ are supported on the discrete set $\mathbb{Z} - \{0\},$ it is enough to fix $k \in \mathbb{Z} - \{0\}$ and show $Ker(A_T)\cap \langle \{\hat{\delta}_{k/T, 1}^{(j)}\}_{j = 1}^\infty \rangle = \{0\}.$ Moreover, because elements of $\mathcal{E}'(\mathcal{H}_\mu)$ are continuous on $C^r(\mathcal{H}_\mu)$ for some $r > 0$, we can fix $r > 0$ and show that any nonzero $\mathcal{D} \in \langle \{\hat{\delta}_{k/T, 1}^{(j)}\}_{j = 1}^r \rangle$ is not in $Ker(A_T)$.  

Fix $k$, $r$ and $\mathcal{D}$.  Then there exists $\{c_j\}_{j = 1}^r$ not all zero such that $$\mathcal{D} = \sum_{j = 1}^r c_j \hat{\delta}_{k/T, 1}^{(j)},$$ and take $\{d_j\}_{j = 1}^{r}$ such that $\sum c_j d_j \neq 0$.    Let $f(\xi) = \frac{e^{2\pi i T \xi} - 1}{2\pi i \xi}$, and we find $h\in C_c^{\infty}(\mathbb{R})$ such that for all $j \geq 1,$ $\delta_{k/T}^{(j)}(f \hat{h}) = d_j.$

Set $\hat{h}(\frac{k}{T}) := \frac{d_1}{f'(\frac{k}{T})},$ and note $f(\frac{k}{T}) = 0,$ so $$\delta_{k/T}^{(1)}(f \hat{h}) = (f'(\frac{k}{T}) \hat{h}(\frac{k}{T}) + f(\frac{k}{T}) \hat{h}'(\frac{k}{T})) = d_1.$$   Notice $f'(\xi) = \frac{2\pi i \xi(2\pi  i T e^{2\pi i T \xi}) - (e^{2\pi i T \xi} - 1) 2\pi i}{(2\pi i \xi)^2},$ so $$f'(\frac{k}{T}) = \frac{(2\pi i k/T) (2\pi i T)}{(2\pi i k/T)^2} = \frac{T^2}{k} \neq 0.$$  Next, for all $2 \leq j \leq r$, we have $$\delta_{k/T}^{(j)}(f \hat{h}) = \sum_{m = 1}^j (_m^j) f^{(m)}(\frac{k}{T}) \hat{h}^{(j - m)}(\frac{k}{T}),$$ which suggests we set $$\hat{h}^{(j - 1)}(\frac{k}{T}) := \frac{d_j - \sum_{m = 2}^{j} (_m^j) f^{(m)}(\frac{k}{T}) \hat{h}^{(j - m)}(\frac{k}{T})}{(_1^j) f'(k/T)}.$$  Then $\delta_k^{(j)}(f \hat{h}) = d_j.$  Now using Taylor series we know a function $\hat{h} \in C_c^{\infty}(\mathbb{R})$ with these derivatives at $\frac{k}{T}$ exists, so $\mathcal{D} \notin Ker(A_T)$.

For $\mu \leq 0$, the Fourier transform is defined along the line $\mathbb{R} + i \subset H$, so for $G\in C^{\infty}(H, d\lambda_\nu)$, define $\hat{G}_1 := \mathcal{F}_1(G(\cdot + i))$.  As before $$(A_T \mathcal{D})(G) = \hat{\mathcal{D}}_1 \int_0^{T} \mathcal{F}_1(G(x + t + i)) dt = \hat{\mathcal{D}}_1 (\frac{e^{2\pi i T \xi} - 1}{2\pi i \xi} \hat{G}_1(\xi)) = 0,$$ whenever $\mathcal{D} \in \langle\{\hat{\delta}_{k/T, 1}\}_{k\in \mathbb{Z}^+}\rangle \left( = \langle\{\hat{\delta}_{k/T, 1}\}_{k\in \mathbb{Z} - \{0\}}\rangle\right)$, so $\langle \{\hat{\delta}_{k/T, 1}\}_{k \in \mathbb{Z}^+} \rangle \subset Ker(A_T).$  

For the other inclusion, let $f, h$ be as in the principal and complementary series cases and define $$G(z) := \int_\mathbb{R} e^{2\pi i (z - i)} \hat{h}(\xi) d\xi.$$  Lemma $\ref{1012}$ shows $G \in C^{\infty}(H, d\lambda_\nu),$ and notice $G(x + i) = h(x),$ which means $\hat{G}_1(\xi) := \mathcal{F}_1(G(\cdot + i)) = \hat{h}(\xi).$  Then we conclude by the same argument.  $ \ \ \Box$\\ \\

Observe $$\int_0^{T} f\circ \phi_t^U dt = A_T f = u\circ \phi_{T}^U - u = \int_0^T \mathcal{L}_U u \circ \phi_t^U dt,$$ which implies $A_T(\mathcal{L}_U u  - f) = 0.$  Lemma $\ref{367}$ now shows the (distributional) $$Ker(A_T) = \bigg\{\begin{array}{rr}\langle \{\hat{\delta}_{k/T}\}_{k \in \mathbb{Z} - \{0\}}\rangle \text{ if } \mu > 0\\ \langle \{\hat{\delta}_{k/T}\}_{k \in \mathbb{Z}^+}\rangle \text{ if } \mu \leq 0.\end{array}$$  So there exists $\{c_k\}_{k \in \mathbb{Z}} \subset \mathbb{C}$ such that \begin{equation}\label{948}\mathcal{L}_U u - f =  \bigg\{\begin{array}{rr}\sum_{k = -\infty}^{\infty} c_k\hat{\delta}_{k/T} \text{ if }\mu > 0\\  \sum_{k = 1}^{\infty} c_k\hat{\delta}_{k/T}\text{ if }\mu \leq 0.\end{array}\end{equation}  

We prove Proposition $\ref{333}$ by showing every $c_k$ is zero in every irreducible representation space.

\begin{lemma}\label{1009}
Let $\mu \geq 1, s > 1$ and $f\in W^{s}(\mathcal{H}_\mu)$.  If there exists $u \in \mathcal{H}_\mu$ such that $$u\circ \phi_T^U - u = A_T f,$$ then $$\mathcal{L}_U u = f.$$
\end{lemma}
$\textit{Proof}:$
By formula ($\ref{948}$), we need to show each coefficient $c_m$ is zero.  Taking Fourier transforms, we have $2\pi i \xi \hat{u} - \hat{f} = \sum_{k = -\infty}^{\infty} c_k \delta_{k/T}.$  Because $u, f \in L^2(\mathbb{R})$, it follows that $\xi\hat{u} - \hat{f} \in L_{loc}^2(\mathbb{R})$.  For all $m \in \mathbb{N}$, let $\{\psi_{j, m}\}_{j\in \mathbb{N}} \subset C_c^{\infty}(\mathbb{R})$ be such that $\psi_{j, m}$ is supported in $[\frac{-1}{j} + \frac{m}{T}, \frac{m}{T} + \frac{1}{j}]$ and $\psi_{j, m}(\frac{m}{T}) = 1 = \| \psi_{j,m}\|_{C^0(\mathbb{R})}$.  Then $$c_m = \lim_{j \rightarrow \infty} \left(\sum_{k = -\infty}^{\infty} c_k \delta_{k/T}\right) (\psi_{j, m}) = \lim_{j \rightarrow \infty} \int_\mathbb{R} (2\pi i\xi\hat{u} - \hat{f})(\xi) \overline{\psi_{j, m}(\xi)} d\xi$$ $$ \leq \lim_{j \rightarrow\infty} \int_{m/T - 1/j}^{m/T + 1/j} |\xi \hat{u}(\xi) - \hat{f}(\xi)| d\xi = 0.$$   Hence, $c_m = 0. \ \ \Box$

For the discrete series case we have the following.

\begin{lemma}\label{1025}
Let $\mu \leq 0, s > 1$ and $f\in W^{s}(H, d\lambda_\nu)$.  If there exists $u \in L^2(H, d\lambda_\nu),$ such that $$u\circ \phi_T^U - u = A_T f,$$ then $$\mathcal{L}_U u = f.$$
\end{lemma}
$\textit{Proof}:$
Let $m \in \mathbb{Z}^+$ and we again show that the coefficient $c_m$ in formula $(\ref{948})$  is zero.  Let $\psi \in C^{\infty}(H, d\lambda_\nu),$ and define $u_y(x) = u(x + i y).$  Observe $\mathcal{L}_U$ is essentially skew-adjoint on $L^2(H, d\lambda_\nu)$, so $$\langle \mathcal{L}_U u, \psi\rangle_{L^2(H, d\lambda_\nu)} = -\langle u, \psi' \rangle_{L^2(H, d\lambda_\nu)}$$ $$ = \int_{0}^{\infty} \left(\int_{-\infty}^\infty u(x + i y) \overline{\psi'(x + i y)} dx \right) y^{\nu - 1} dy$$ \begin{equation}\label{1015} = \int_0^{\infty} \langle \hat{u}_y(\xi), (-2\pi i \xi) \hat{\psi}_y\rangle_{L^2(\mathbb{R})} y^{\nu - 1} dy = \langle (2\pi i \xi) \mathcal{F}_1 u, \mathcal{F}_1 \psi\rangle.\end{equation} 

Let $\epsilon > 0$ and let $\hat{g}_{\epsilon, m} \in C_c^{\infty}(\mathbb{R})$ be a bump function supported on $[\frac{m}{T} - \epsilon, \frac{m}{T} + \epsilon]$ that satisfies $\|\hat{g}\|_{C^0(\mathbb{R})} = e^{-2\pi m/T}$ and $\hat{g}_{\epsilon, m}(\frac{m}{T}) = e^{-2\pi m/T}$.  Now define \begin{equation}\label{1017}G_{\epsilon, m}(z) := \int_\mathbb{R} e^{2\pi i (z - i) \xi} \hat{g}_{\epsilon, m}(\xi) d\xi = \int_\mathbb{R} e^{2\pi i x \xi} e^{-2 \pi (y - 1) \xi} \hat{g}_{\epsilon, m}(\xi) d\xi\end{equation} and note that $G_{\epsilon, m}(x + i) = g_{\epsilon, m}(x).$  So for $k\in \mathbb{Z}^+$, $$\hat{\delta}_{k/T}(G_{\epsilon, m}) = \int_\mathbb{R} e^{-2\pi i k/T (x + i)} G_{\epsilon, m}(x + i) dx$$ $$ = e^{2\pi k/T} \int_\mathbb{R} e^{2\pi i x k/T} g_{\epsilon, m}(x) dx = e^{2\pi k/T} \hat{g}_{\epsilon, m}(\frac{k}{T}) = \bigg\{\begin{array}{rr} 1\text{ if }k = m \\ 0\text{ if }k \neq m.\end{array}$$  

So for any $\epsilon > 0$, $$|c_m| = |\left(\sum_{k = 1}^{\infty} c_k \hat{\delta}_{k/T}\right) G_{\epsilon, m}| \leq |\langle \mathcal{L}_U u, G_{\epsilon, m}\rangle_{L^2(H, d\lambda_\nu)}| + |\langle f, G_{\epsilon, m}\rangle_{L^2(H, d\lambda_\nu)}|.$$  Because $f \in L^2(H, d\lambda_\nu)$, note $\lim_{\epsilon \rightarrow 0} |\langle f, G_{\epsilon, m}\rangle_{L^2(H, d\lambda_\nu)}| = 0,$ so it suffices to prove $\lim_{\epsilon\rightarrow 0}|\langle \mathcal{L}_U u, G_{\epsilon, m}\rangle_{L^2(H, d\lambda_\nu)}| = 0.$

Set $u_y(x) = u(x + i y)$ and let $G_{\epsilon, m, y}$ be defined similarly.  Using $(\ref{1015})$, we have $$|\langle \mathcal{L}_U u, G_{\epsilon, m}\rangle_{L^2(H, d\lambda_\nu)}| = |\langle (2\pi i \xi) \mathcal{F}_1 u, \mathcal{F}_1 G_{\epsilon, m}\rangle_{L^2(H, d\lambda_\nu)}|$$ \begin{equation}\label{4011} = |2\pi \int_0^{\infty} \langle \hat{u}_y, \xi \hat{G}_{\epsilon, m, y}\rangle_{L^2(\mathbb{R})} y^{\nu - 1} dy| \leq 2\pi \int_0^{\infty} \| \hat{u}_y\|_{L^2(\mathbb{R})} \|\xi \hat{G}_{\epsilon, m, y}\|_{L^2(\mathbb{R})} y^{\nu - 1} dy.\end{equation}  A calculation now proves $\|\xi \hat{G}_{\epsilon, m, y}\|_{L^2(\mathbb{R})} \leq C_m \sqrt{\epsilon} e^{-2\pi(y - 1)(m/T - \epsilon)}.$  As $m \geq 1$, we conclude $c_m = \lim_{\epsilon\rightarrow 0} (\ref{4011}) = 0. \ \ \Box$ 

Lastly, we handle the complementary series.

\begin{lemma}\label{1013}
Let $0 < \mu < 1, s > 1$ and $f\in W^{s}(\mathcal{H}_\mu)$.  If there exists $u \in \mathcal{H}_\mu,$ such that $u\circ \phi_T^U - u = A_T f,$ then $\mathcal{L}_U u = f.$
\end{lemma}
$\textit{Proof}:$
Let $\hat{g} \in C_c^{\infty}(\mathbb{R})$ be a bump function supported on $[-1, 1]$ with $\hat{g}(0) = 1.$  Additionally, fix $m \in \mathbb{Z}$, and for all $n \in \mathbb{Z}^+$, define $\hat{g}_{m, n}(\xi) := \hat{g}(n(\xi - m)).$  Notice that $g_{m, n}(x) = \frac{1}{n} e^{2\pi i m x} g(\frac{x}{n}).$  Then from formula $(\ref{948})$, $$|c_m| = |\left(\sum_{k = -\infty}^{\infty} c_k \hat{\delta}_{k/T}\right)(g_{m, n})| \leq |\langle \mathcal{L}_U u, g_{m, n}\rangle_{\mathcal{H}_\mu}| + |\langle f, g_{m, n}\rangle_{\mathcal{H}_\mu}|$$ \begin{equation}\label{943} \leq |\langle u, g_{m, n}'\rangle_{\mathcal{H}_\mu}| + |\langle f, g_{m, n}\rangle_{\mathcal{H}_\mu}| \leq \| u \|_{\mathcal{H}_\mu} \| g_{m, n}' \|_{\mathcal{H}_\mu} + \| f \|_{\mathcal{H}_\mu} \| g_{m, n}\|_{\mathcal{H}_\mu}.\end{equation}

We will estimate $\| g_{m, n}' \|_{\mathcal{H}_\mu}$ and $\| g_{m, n} \|_{\mathcal{H}_\mu}$ using the following lemma.

\begin{lemma}\label{1095}
Let $0 < \mu < 1$ and $h \in \mathcal{H}_\mu$.  Then for any $q \in [1, \frac{1}{\nu})$, there exists a constant $C_{q, \nu} > 0$ such that $$\| h \|_{\mathcal{H}_\mu}^2 \leq C_{q, \nu} \| h \|_{L^1(\mathbb{R})} (\| h \|_{L^{q}(\{|x| \geq 1\})} + \| h \|_{L^{\infty}(\{|x| < 1\})}).$$
\end{lemma}
$\textit{Proof}:$
For a given function $h \in \mathcal{H}_\mu$, \begin{equation}\label{942}\| h \|_{\mathcal{H}_\mu}^2 = \int_{\mathbb{R}^2} \frac{h(r + x) h(x)}{|r|^{1 - \nu}} dr dx = \int_{\mathbb{R}} h(x) \left(\int_\mathbb{R} \frac{h(r + x)}{|r|^{1 - \nu}} dr\right) dx.\end{equation}  Observe that $$\int_\mathbb{R} \frac{h(r + x)}{|r|^{1 - \nu}} dr \leq \int_{\{|r| \geq 1\}} \frac{|h(r + x)|}{|r|^{1 - \nu}} dr + \int_{\{|r| < 1\}} \frac{|h(r + x)|}{|r|^{1 - \nu}} dr,$$ and notice \begin{equation}\label{945}\int_{\{|r| < 1\}} \frac{|h(r + x)|}{|r|^{1 - \nu}} dr \leq C_\nu \| h \|_{L^\infty(\mathbb{R})}.\end{equation}  Additionally, let $p, q > 0$ satisfy $\frac{1}{p} + \frac{1}{q} = 1.$  Then Holder's inequality gives $$\int_{\{|r| \geq 1\}} \frac{|h(r + x)|}{|r|^{1 - \nu}} dr \leq \left(\int_{\{|r| \geq 1\}} \frac{1}{|r|^{p(1 - \nu)}} dr\right)^{1/p} \left(\int_{\{|r| \geq 1\}} |h(r + x|)|^{q} dr\right)^{1/q}.$$  Let $1 < p \leq \infty$ be such that $p(1 - \nu) > 1.$  So $(1 - \nu) > \frac{1}{p} = 1 - \frac{1}{q},$ and therefore 
$1 \leq q < \frac{1}{\nu}.$  

When $p = \infty$, then $q = 1$ and $$(\ref{942}) \leq C_{\nu} \| h \|_{L^1(\mathbb{R})} (\| h \|_{L^{1}(\{|x| \geq 1\})} + \| h \|_{L^{\infty}(\{|x| < 1\})}).$$  When $p < \infty$, then $1 < q < \frac{1}{\nu}$ and we conclude \begin{equation}\label{946}\int_{\{|r| \geq 1\}} \frac{|h(r + x)|}{|r|^{1 - \nu}} dr \leq C_{q, \nu} \|h \|_{L^{q}(\mathbb{R})}.\end{equation}  
Then $$(\ref{942}) \leq C_{q, \nu} \| h \|_{L^1(\mathbb{R})} (\| h \|_{L^{q}(\{|x| \geq 1\})} + \| h \|_{L^{\infty}(\{|x| < 1\})}). \ \ \Box$$  

We now know there exists $q > 1$ such that $$(\ref{943}) \leq C_{\nu, q}\| u \|_{\mathcal{H}_\mu} \| g_{m, n}'\|_{L^1} (\|g_{m, n}' \|_{L^{q}} + \|g_{m, n}' \|_{L^\infty})$$ \begin{equation}\label{962} + \|f \|_{\mathcal{H}_\mu} \| g_{m, n}\|_{L^1} (\|g_{m, n} \|_{L^{q}} + \|g_{m, n} \|_{L^\infty}).\end{equation}

Notice $$\frac{d}{dx} g_{m, n}(x) = \frac{1}{n^2} g_{m, n}'(\frac{x}{n}) e^{2\pi i m x} + \frac{2\pi i m}{n} e^{2\pi i m x} g(\frac{x}{n}).$$   Then there is a constant $C > 0$ such that $$\| g_{m, n}'\|_{L^1} + \|g_{m, n}\|_{L^1} \leq C m, \ \|g_{m, n}'\|_{L^{q_1}} +  \|g_{m, n}\|_{L^{q_1}} \leq C m n^{1/q - 1},$$ $$ \ \|g_{m, n}' \|_{L^\infty} + \|g_{m, n} \|_{L^\infty}  \leq C \frac{m}{n}.$$  Therefore, $\lim_{n\rightarrow\infty} (\ref{962}) = 0,$ and we conclude $c_m = 0. \ \ \Box$ \\

$\textit{Proof of Proposition }\ref{333}$:  This follows from Lemmas $\ref{1009}, \ref{1025}$ and $\ref{1013}. \ \ \Box$

\section{ }
\subsection{Proof of Lemma $\ref{1080}$}
The first step in proving Lemma $\ref{1080}$ is the following.

 \begin{lemma}\label{715}
Let $\mu \in spec(\Box)$. There exists a dual set of functions $\{f_n\} \subset C^{\infty}(\mathcal{K}_\mu)$ to the spanning set of distributions $\{\mathcal{D}_n\}\subset \mathcal{I}(\mathcal{K}_\mu)$ so that for all $n$, $$\mathcal{D}_n(f_\tau) = \Big\{\begin{array}{rr}1 \text{ if }\mathcal{D}_n = \mathcal{D}_\tau \\0 \text{ otherwise. }\end{array}$$ 

Moreover, for $s > 1$, there is a constant $C_s > 0$ such that for all $\tau\in \mathbb{Z}$, $$|c_\tau(x, N, 1, s)| \leq C_s (1 + |\nu|)^s |\tau|^{s}.$$
\end{lemma}

This will be immediate from the next two lemmas.
 \begin{lemma}\label{1060}
Lemma $\ref{715}$ holds when $\mu \in spec(\triangle_{SM})$.
\end{lemma}
$\textit{Proof}:$
First suppose that $s, \tau \in \mathbb{Z}^+$.  Let $\hat{f}_0 \in C_c^{\infty}(\mathbb{R})$ be supported in $(-1/2, 1/2)$, and let $\hat{f}_0(0) = 1$.  Now let $\chi_\tau = e^{2\pi i \tau x} f_0 \in C_c^{\infty}(\mathbb{R})$, so that \begin{equation}\label{2102}\hat{\delta}_k(\chi_\tau) = \Big\{\begin{array}{rr} 1 \text{ if }k = \tau\\ 0 \text{ otherwise} \end{array}.\end{equation}  Note that because $\chi_\tau \in C_c^{\infty}(\mathbb{R})$, we know that $\delta^{(0)}(\chi_\tau) = 0.$  Now let $\pi: W^1(SM)\rightarrow W^1(SM)\cap (Ann(\mathcal{I}(SM)))^\bot$ be orthogonal projection, so $\chi_\tau = \pi\chi_\tau \oplus (1 - \pi)\chi_\tau$ and $(\ref{2102})$ still holds for $\pi \chi_\tau$, and $\delta^{(0)}(\pi \chi_\tau) = 0.$ 

For the estimate, Sobolev embedding gives $$|c_\tau| = |c_\tau \mathcal{D}_\tau(Q_\mu^{-1} \pi\chi_\tau)|$$ $$ \leq |\left(\left(\sum_{n = -\infty}^{\infty} c_n \mathcal{D}_n + d_0 \mathcal{D}^{0}\right)\oplus \mathcal{R}\right)(Q_\mu^{-1} \pi \chi_\tau)| + |\mathcal{R}(Q_\mu^{-1} \pi\chi_\tau)|$$ \begin{equation}\label{4012} = \frac{1}{N}\sum_{k = 0}^{N - 1} Q_\mu^{-1} \pi \chi_n(\phi_n^U(x_0)) + 0 \leq C_s \| \chi_\tau \|_s.\end{equation} We can estimate $\|\chi_\tau \|_s$ using our concrete formulas for $X, Y, \Theta$ in $(1 + \triangle)^{s/2} = (1 - (X^2 + Y^2 + \Theta^2))^{s/2}.$  Because $supp(f_0) \subset (-1/2, 1/2)$, we have \begin{equation}\label{2101}\| \chi_\tau \|_s \leq C_s (1 + |\nu|)^s \| (1 + U^s) \chi_\tau \|_{\mathcal{H}_\mu} \leq 
C_s (1 + |\nu|)^s (1 + |\tau|)^s \| f_0\|_s.\end{equation}

Then the lemma for real $s > 1$ follows by interpolation \cite{11}.  $\ \ \ \Box$
\begin{lemma}\label{1012}
Lemma $\ref{715}$ holds when $\mu \leq 0$.  
\end{lemma}
$\textit{Proof}:$
By Lemma $\ref{767}$, there is nothing to prove if $\tau \leq 0$, so let $\tau \in \mathbb{Z}^+$.  First say $s \in \mathbb{Z}^+$.  Let $\hat{g}_\tau \in C_c^{\infty}(\mathbb{R})$ be supported in $[-\frac{1}{2} + \tau, \tau + \frac{1}{2}]$ and satisfy $\hat{g}_\tau(\tau) = \frac{1}{e^{2\pi \tau}}.$  Define $G_\tau(z)$ as in $(\ref{1017})$ so that $G_\tau(x + i y) = \mathcal{F}_1^{-1}h(\xi, y)$, where $h(\xi, y) = e^{-2\pi (y - 1)\xi} g_\tau(\xi),$ and $G_\tau$ is holomorphic.  Then integration by parts proves $\delta^{(0)}(G_\tau) = 0$, so by $(\ref{4012})$, we only need an estimate on $\|G_\tau\|_s.$  For this, integration by parts proves \begin{equation}\label{2100}\|G_\tau\|_s \leq C_s (1 + |\nu|)^s \tau^s\end{equation} when $s \in \mathbb{Z}^+$, and the estimate for any $s > 1$ follows by interpolation \cite{11}. $ \ \ \Box$

$\textit{Proof of Lemma }\ref{715}:$ This is immediate from Lemmas $\ref{1060}$ and $\ref{1012}. \ \ \Box$

$\textit{Proof of Lemma }\ref{1080}:$
We prove the identity for $c_\tau(x_0, N, 1, s)$.  
Notice \begin{equation}\label{1087}\int_{0}^1e^{2\pi i \tau t} \phi_{-t}^U\left(\left(\sum_{k = -\infty}^{\infty} c_k(x_0, N, 1, s) \mathcal{D}_k + d_0(x_0, N, 1, s)\mathcal{D}^{0}\right) \oplus \mathcal{R}(x_0, N, 1, s)\right) dt\end{equation} $$ = \int_0^1e^{2\pi i \tau t} \phi_{-t}^U\left(\frac{1}{N}\sum_{n = 0}^{N - 1} \delta_{\phi_n^U(x_0)} \right) dt = \frac{1}{N}\int_0^N e^{2\pi i \tau t} (\phi_t^U (x_0))^* dt.$$  Using Proposition $\ref{156}$, we get $$\int_0^1 e^{2\pi i \tau t} \phi_{-t}^U\left(d_0(x_0, N, 1, s)\mathcal{D}^{0} + \mathcal{R}(x_0, N, 1, s)\right) dt\in W^{-s}(SM),$$ and Sobolev embedding shows $$\frac{1}{N}\int_0^N e^{2\pi i \tau t} (\phi_t^U(x_0))^* dt\in W^{-s}(SM).$$  Hence, $$\int_{0}^1e^{2\pi i \tau t} \phi_{-t}^U\left(\sum_{k = -\infty}^{\infty} c_k(x_0, N, 1, s) \mathcal{D}_k\right) dt\in W^{-s}(SM).$$  Then
 we may separate the integral $(\ref{1087})$ and conclude $$\int_0^1 e^{2\pi i \tau t} c_{\tau}(x_0, N, 1, s) \phi_{-t}^U\mathcal{D}_0 dt = \frac{1}{N}\int_0^N e^{2\pi i \tau t}(\phi_t^U(x_0))^* dt$$ \begin{equation}\label{1050} - \int_0^1e^{2\pi i \tau t} \phi_{-t}^U\left(\sum_{n \neq \tau} c_n \mathcal{D}_n\right) dt - d_0(x_0, N, 1, s)  \int_0^1e^{2\pi i \tau t} \phi_{-t}^U \mathcal{D}^0 dt - \int_0^1 e^{2\pi i \tau t} \phi_{-t}^U \mathcal{R}dt.\end{equation}

In the same way one shows $\sum_{k  \neq \tau} c_k \mathcal{D}_k \in W^{-s}(\mathcal{K}_\mu).$  Let $f\in C^{\infty}(\mathcal{K}_\mu)$.  Then because the unitary equivalence $Q_\mu$ intertwines $\phi_t^U$, we get $$(\phi_{-t}^U\mathcal{D}_k)(f) = (Q_\mu^* \hat{\delta}_{k})(f \circ \phi_{t}^U)  = \hat{\delta}_{k}(Q_\mu(f\circ\phi_{t}^U )) = \hat{\delta}_k((Q_\mu f) \circ \phi_{t}^U)
$$ \begin{equation}\label{163} = \hat{\delta}_{k}(Q_\mu f (\cdot - t)) = e^{-2\pi i k t} \hat{\delta}_{k}(Q_\mu f) = e^{-2\pi i k t} \mathcal{D}_k(f).\end{equation}  Lemma $\ref{715}$ shows that $|c_k| \leq C_s (1 + |\nu|)^s (1 + |k|)^{s}$ for all $k$ and Lemma $\ref{787}$ proves the Fourier transform $\mathcal{F}(Q_\mu f)$ decays faster than the reciprocal of any polynomial, so we get $$\sum_{k \in \mathbb{Z}} c_k \mathcal{D}_k(f) = \sum_{k \in \mathbb{Z}} c_k \hat{\delta}_k(Q_\mu f)$$ converges absolutely.  Then using $(\ref{163})$, \begin{equation}\label{169}\int_0^1e^{2\pi i \tau t} \phi_{-t}^U \sum_{k \neq \tau} c_k \mathcal{D}_k (f)dt = \sum_{k \neq \tau} c_k \int_0^1e^{2\pi i \tau t} \phi_{-t}^U \mathcal{D}_k (f) dt = 0.\end{equation}    

Moreover, because $\mathcal{D}^0$ is flow invariant, $$d_0(x_0, N, 1, s)  \int_0^1e^{2\pi i \tau t} \phi_{-t}^U \mathcal{D}^0 dt = 0.$$  Then combining this with $(\ref{1050})$, $(\ref{163})$ and $(\ref{169})$ gives the identity for $c_\tau$.  The arguments for the coefficients $c_0$ and $d_0$ are the same.  $\ \ \Box$

\subsection{Proof of Lemma $\ref{787}$}
$\textit{Proof of Lemma }\ref{787}:$  This will follow from the next two lemmas. 
\begin{lemma}\label{1040}
Let $\mu \in spec(\triangle_M), r \geq 2$ and $f\in W^{3r + 2}(\mathcal{H}_\mu).$  Then for all $\xi\in \mathbb{R}$, $$|\hat{f}(\xi)| \leq C_r \|f\|_{W^{3r + 2}(\mathcal{H}_\mu)} (1 + |\xi|)^{-r}.$$  
\end{lemma}
$\textit{Proof}:$
Observe that $$\| \xi^r \hat{f} \|_{C^0(\mathbb{R})} \leq C (\| \xi^r \hat{f} \|_{L^2(\mathbb{R})} + \|\frac{d}{d\xi} (\xi^r \hat{f}) \|_{L^2(\mathbb{R})})$$ \begin{equation}\label{786}\leq C_r (\| f^{(r -1)}\|_{L^2(\mathbb{R})} + \| f^{(r)} \|_{L^2(\mathbb{R})} + \|x f^{(r)} \|_{L^2(\mathbb{R})}).\end{equation}  By the change of variable $x = \tan(\theta)$, recall $U = -\cos^2(\theta)\frac{d}{d\theta}.$  Then by induction on integers $r \geq 0$, we have $$|f^{(r)}(x)| = |U^r f(\theta)| \leq C_r \sum_{j = 1}^r \cos^{r + j}(\theta) |f^{(j)}(\theta)|.$$  We abuse notation slightly and write $f(x) = f(\theta) = \Phi(\theta)\cos^{1 + \nu}(\theta),$ where $\Phi(\theta) = \sum_{m = -\infty}^{\infty} c_m e^{2\pi i m \theta}$.  One likewise proves $$|f^{(j)}(\theta)| \leq C_r (1 + |\nu|)^j \sum_{k = 0}^j |\Phi^{(k)}(\theta)| \cos^{1 - j + \nu}(\theta).$$  Now Sobolev's inequality implies $|\Phi^{(k)}(\theta)| \leq $ $\| f \|_{W^{k + 1}(\mathcal{H}_\mu)}$.  Combining and switching to $\mathbb{R}$ coordinates, we find a constant $C_r > 0$ such that for all $x\in \mathbb{R}$, $$|f^{(r)}(x)| \leq C_r (1 + |\nu|)^r (1 + |x|)^{-(r + 1 + \nu)} \| f \|_{r + 1} \leq C_r (1 + |x|)^{-(r + 1 + \nu)} \| f \|_{2r + 1},$$ where the second inequality follows from Lemma 6.3 of \cite{8}.  Therefore, $(\ref{786}) \leq C_r \| f \|_{2r + 1}. \ \ \Box$

\begin{lemma}\label{1037}
Let $\mu \leq 0,$ $r \geq 2$ and $f\in W^{3r + 2}(H, d\lambda_\nu)$.  Then there is a constant $C_r > 0$ such that $$|\mathcal{F}_1 f(\xi + i)| \leq C_r \| f \|_{3r + 2} (1 + |\xi|)^{-r}.$$
\end{lemma}
$\textit{Proof}:$
This time using formula $(\ref{350})$ and Lemma $\ref{447}$, we conclude. $\ \ \Box$

Now Lemma $\ref{787}$ is immediate from Lemmas $\ref{1040}$ and $\ref{1037}. \ \ \Box$

\end{document}